\documentclass{amsart}
\usepackage{amssymb,latexsym}

\numberwithin{equation}{section}

\newtheorem{theorem}{Theorem}[section]
\newtheorem{proposition}[theorem]{Proposition}
\newtheorem*{proposition*}{Proposition}
\newtheorem{lemma}[theorem]{Lemma}
\newtheorem*{lemma*}{Lemma}
\newtheorem{corollary}[theorem]{Corollary}
\newtheorem*{corollary*}{Corollary}
\newtheorem{fact}[theorem]{Fact}
\newtheorem*{fact*}{Fact}

\theoremstyle{definition}
\newtheorem{definition}[theorem]{Definition}
\newtheorem*{notterm}{Notation--Terminology}

\newtheorem*{example*}{Example}

\theoremstyle{remark}

\newtheorem*{remark*}{Remark}

\newtheorem*{note*}{Note}

\newtheorem*{notation}{Notation}

\newcommand{\hf}{\widehat{f}}
\newcommand{\hg}{\widehat{g}}
\newcommand{\hh}{\widehat{h}}
\newcommand{\hmu}{\widehat{\mu}}
\newcommand{\hnu}{\widehat{\nu}}
\newcommand{\hA}{\widehat{A}}

\newcommand{\hG}{\widehat{G}}

\newcommand{\hK}{\widehat{K}}
\newcommand{\hF}{\widehat{F}}
\newcommand{\hC}{\widehat{C}}
\newcommand{\hS}{\widehat{S}}

\newcommand{\hW}{\widehat{W}}
\newcommand{\tkappa}{\tilde{\kappa}}

\newcommand{\abs}[1]{\lvert#1\rvert}
\newcommand{\bigabs}[1]{\bigl\lvert#1\bigr\rvert}
\newcommand{\Bigabs}[1]{\Bigl\lvert#1\Bigr\rvert}
\newcommand{\biggabs}[1]{\biggl\lvert#1\biggr\rvert}

\newcommand{\set}[1]{\lbrace#1\rbrace}
\newcommand{\bigset}[1]{\bigl\lbrace#1\bigr\rbrace}

\newcommand{\biggset}[1]{\biggl\lbrace#1\biggr\rbrace}
\newcommand{\Biggset}[1]{\Biggl\lbrace#1\Biggr\rbrace}
\newcommand{\norm}[1]{\lVert#1\rVert}
\newcommand{\bignorm}[1]{\bigl\lVert#1\bigr\rVert}
\newcommand{\Bignorm}[1]{\Bigl\lVert#1\Bigr\rVert}

\newcommand{\Biggnorm}[1]{\Biggl\lVert#1\Biggr\rVert}
\newcommand{\ls}{\leqslant}
\newcommand{\gr}{\geqslant}

\def \1{\textup{(}\textit{i\/}\textup{)}}
\def \2{\textup{(}\textit{ii\/}\textup{)}}
\def \3{\textup{(}\textit{iii\/}\textup{)}}
\def \4{\textup{(}\textit{iv\/}\textup{)}}

\begin{document}

\title{On Mixing and Ergodicity in Locally Compact Motion Groups}

\author{M.\ Anoussis}
\address{M.\ Anoussis: Department of Mathematics, University of the Aegean,
832 00 Karlovasi --- Samos, Greece}
\email{mano@aegean.gr}

\author{D.\ Gatzouras}
\address{D.\ Gatzouras: Laboratory of  Mathematics and Statistics,
Agricultural University of Athens, Iera Odos 75, 118 55 Athens,
Greece}
\email{gatzoura@aua.gr}

\subjclass[2000]{Primary; Secondary}
\date{\today}

\begin{abstract}
Let $G$ be a semi-direct product $G=A\times_\varphi K$ with $A$
Abelian and $K$ compact.  We characterize spread-out probability
measures on $G$ that are mixing by convolutions by means of their
Fourier transforms. A key tool is a spectral radius formula for
the Fourier transform of a regular Borel measure on $G$ that we
develop, and which is analogous to the well-known
Beurling--Gelfand spectral radius formula.  For spread-out probability measures on $G$, we
also characterize ergodicity by means of the Fourier transform of
the measure. Finally, we show that spread-out probability measures
on such groups are mixing if and only if they are weakly mixing.
\end{abstract}

\maketitle

\section{Introduction}
\label{intro}

The purpose of this paper, which may be viewed as a sequel to
\cite{AG}, is to exploit methods of non-commutative harmonic
analysis to study random walks on locally compact groups.

Our starting point is the following spectral radius formula for a
regular (complex) Borel measure $\mu$ on a locally compact
Hausdorff group $G$:
\begin{equation}\label{SRFM}
\lim\limits_{n\to\infty}\norm{\mu^n}^{1/n} =
\sup\limits_{U}\varrho\bigl(\hmu(U)\bigr) \vee
\inf\limits_{n\in\mathbb{N}}
\norm{(\mu^n)_{\mathrm{s}}}^{1/n},\end{equation} where
$\mu^n:=\mu\ast\cdots\ast\mu$ denotes $n$-fold convolution of
$\mu$ with itself, $(\mu^n)_{\mathrm{s}}$ is the singular part of
$\mu^n$ with respect to Haar measure $\lambda_G$ and $U$ runs
through a complete set of continuous irreducible unitary
representations of $G$ (also $a\vee b:=\max\set{a,b}$); here, and
throughout the paper, $\lambda_G$ is a fixed left Haar measure on
$G$, $\hmu(U)$ denotes the Fourier transform of $\mu$ at the
unitary representation $U$, and $\varrho\bigl(\hmu(U)\bigr)$ its
spectral radius.  When $G$ is Abelian, this formula is a direct
consequence of Gelfand theory for the commutative Banach algebra
$M(G)$ of regular (complex) Borel measures on $G$, and for groups
$G$ with a symmetric group algebra $L^1(G)$, \eqref{SRFM} has been
established by Palmer \cite{PaA} for absolutely continuous
measures $\mu$.

In \cite{AG} \eqref{SRFM} was established for compact (Hausdorff)
groups $G$ and then used to study random walks on such groups. In
particular, one of its uses there was to characterize those
regular Borel probability measures $\mu$ on the compact group $G$,
for which $\mu^n\to\lambda_G$ in the total variation norm (here we
assume that $\lambda_G$ has been chosen to satisfy
$\lambda_G(G)=1$). Of course when $G$ is non-compact, $\mu^n$
cannot converge to Haar measure for any probability measure $\mu$,
and neither, of course, can $(1/n)\sum_{k=0}^{n-1}\mu^k$. Two
conditions that have been studied extensively instead are
\begin{equation}\label{MC}
\lim\limits_{n\to\infty}\norm{f\ast\mu^n}_1=0\qquad \forall\, f\in
L^1_0(G)\end{equation} and
\begin{equation}\label{EC}
\lim\limits_{n\to\infty}\Biggnorm{\frac{1}{n}\sum\limits_{l=0}^{n-1}
f\ast\mu^k}_1=0\qquad \forall\, f\in L^1_0(G),\end{equation} where
$L^1_0(G)$ denotes the closed two-sided ideal of functions $f\in
L^1(G)$ with $\int_G f\,d\lambda_G=0$.  We shall call probability
measures $\mu\in M(G)$ satisfying \eqref{MC} \emph{mixing by
convolutions\/}, and those satisfying \eqref{EC} \emph{ergodic by
convolutions\/}, adhering to terminology introduced by Rosenblatt
in \cite{Rblatt}, in view of the fact that $\mu$ satisfies
\eqref{MC} if and only if the associated random walk is mixing,
and similarly for ergodicity. The results of this paper are then
as follows:

\begin{itemize}
\item[(A)]  We establish the spectral radius formula \eqref{SRFM}
for arbitrary regular Borel measures in motion groups (Theorem
\ref{srfm}).
\item[(B)]  Using (A), we show that in a motion group $G=A\times_\varphi K$ with $G$
acting regularly on $\hA$, a spread-out probability measure
$\mu\in M(G)$ is mixing by convolutions if and only if
\begin{equation}\label{SR}\varrho\bigl(\hmu(U)\bigr)< 1\qquad
\forall\, [U]\in\hG\smallsetminus\set{\mathbf{1}_G},\end{equation}
where $\hG$ is the unitary dual of $G$ and $\mathbf{1}_G$
designates the trivial representation of $G$ (Theorem \ref{main}).
\item[(C)] Under the same conditions on $G$, we show that a spread-out probability measure
$\mu\in M(G)$, is ergodic by convolutions if and only if
\begin{equation}\label{Sp} 1\notin\sigma\bigl(\hmu(U)\bigr)\qquad
\forall\, [U]\in\hG\smallsetminus\set{\mathbf{1}_G},\end{equation}
where $\sigma\bigl(\hmu(U)\bigr)$ denotes the spectrum of the
operator $\hmu(U)$ (Theorem \ref{main-erg}).
\end{itemize}
Finally, as a consequence of our approach, we are able to address
a query in \cite{Rblatt} (p.~33).  We show that in a motion group
$G=A\times_\varphi K$ with $G$ acting regularly on $\hA$, a
spread-out probability measure $\mu\in M(G)$ is mixing by
convolutions if and only if it satisfies the apparently weaker
condition
\begin{equation}\label{WMC}
\frac{1}{n}\sum\limits_{k=0}^{n-1} \bigabs{\int_G (f \ast \mu^k)
h\, d\lambda_G}\to 0\qquad \forall\,f\in
L^1_0(G)\quad\forall\,h\in L^\infty(G)\end{equation} (Corollary
\ref{mix=wmix}).  Adhering to standard terminology again, we shall
call probability measures satisfying \eqref{WMC} \emph{weakly
mixing by convolutions\/}.  Besides Abelian and compact groups,
the only other groups we know of for which mixing has been shown
to be equivalent to weak-mixing are groups possessing small
invariant neighborhoods (SIN); this is done for arbitrary regular
Borel probability measures in a recent paper by Jaworski
(\cite{Jaw2}).

Let us expand briefly and also motivate these results. By
\emph{motion group\/} we shall mean a group which is a semi-direct
product $G=A\times_\varphi K$ with $A$ Abelian and $K$ compact;
both $A$ and $K$ are assumed to be locally compact and Hausdorff
here. In all our results concerning motion groups we shall also
assume that $G$ acts regularly on the dual group $\hA$ of $A$
\cite[p.~183]{Fol}; this condition is automatically satisfied when
$G$ is second countable (\cite{Glimm}). Recall also that a
probability measure $\mu\in M(G)$ on a locally compact Hausdorff
group is called \emph{spread-out\/} if not all of its convolution
powers $\mu^n$ are singular with respect to Haar measure
$\lambda_G$.

When $G$ is an Abelian group, \eqref{SR} reduces to
$\abs{\hmu(\chi)}<1$ for all non-trivial characters $\chi$ of $G$,
and (B) in this case is a result of Foguel \cite{Fog}; also, for
Abelian $G$, \eqref{Sp} simply says that $\hmu(\chi)\neq 1$ for
all non-trivial characters $\chi$ of $G$, and (C) in this case
follows from the Choquet--Deny theorem \cite{ChD} (see also
\cite{Ram} and \cite{Rblatt}).  On the other hand, when $G$ is
non-commutative, the Fourier transform $\hmu(U)$ of a probability
measure $\mu$ is an operator on a Hilbert space, so it is not
immediately clear what the appropriate generalizations of these
conditions are, to begin with.  A natural choice is to try to use
the norms of the operators $\hmu(U)$ to give conditions for mixing
and ergodicity (see \cite{JRT}, \cite{Kan}), but as it turns out,
$\norm{\hmu(U)}$ does not characterize mixing nor ergodicity.  Our
proof of (B) uses the spectral radius formula \eqref{SRFM},
through which the spectral radius $\varrho\bigl(\hmu(U)\bigr)$
emerges naturally.  Let us also remark that (B), (C), and the
equivalence of mixing and weak mixing, hold without the spread-out
assumption on $\mu$ when $G$ is either Abelian or compact, and
that the results for $G$ compact, although not explicitly
appearing in the existing literature, once appropriately
formulated, also follow from existing results, notably the work of
Kawada and Ito \cite{Kaw-Ito} (see Section \ref{cr}).

Our result (C) on ergodicity is closely related to a result of
Jaworski \cite{Jaw1}, who shows that in a locally compact second
countable group of polynomial growth, a spread-out probability
measure is ergodic by convolutions if and only if it is adapted
(see Section \ref{cr}).  In fact, in Section \ref{cr}, we give a
short direct argument showing how our Theorem \ref{main-erg}, for
second countable motion groups, may be obtained from Jarowski's
result, whose proof relies on structure theory for groups of
polynomial growth.  As adaptedness of a probability measure is
known to not be equivalent to ergodicity in general groups however
(see Rosenblatt \cite{Rblatt}), condition \eqref{Sp} may well be
worth considering, especially in view of the discussion of
Rosenblatt's example in Section \ref{cr}. For the same reason, we
have retained a proof of Theorem \ref{main-erg} relying solely on
the methods of the present paper.

Finally, let us also mention that we in fact obtain stronger
results in one direction in (B) and (C): \eqref{SR} is necessary
for weak mixing, and hence also for mixing, and \eqref{Sp}
necessary for ergodicity, for spread-out probability measures in
\emph{any\/} CCR group (Proposition \ref{mix=wmix} and Corollary
\ref{cor-erg}; see also Corollary \ref{cor1}).

\bigskip

\noindent In \cite{Ram} Ramsey and Weit give different proofs of
Foguel's result on mixing and the Choquet--Deny theorem for
Abelian groups, which are more illuminating from the point of view
of harmonic analysis.  We now briefly contrast the proof of the
more involved direction of our result on mixing, namely the
sufficiency of the condition \eqref{SR} in (B) above, to the
corresponding proof of Ramsey and Weit for Abelian groups $G$.
This will also indicate how the spectral radius formula
\eqref{SRFM} for measures, rather than functions in $L^1$, is
relevant. The proof of Ramsey and Weit relies on the fact that if
$f\in L^1(G)$ is such that $\hf$ has compact support not
containing $0$, then $f$ factorizes as $f\ast h=f$, with $h\in
L^1(G)$ and such that $\hh$ has compact support not containing $0$
. They then use the Beurling--Gelfand spectral radius formula for
functions in $L^1$ and the fact that $\set{\mathbf{1}_G}$  is a
set of synthesis to conclude the proof.

In our setting, we work with the representations
$\varLambda_\alpha,\ \alpha\in\hA$, where $\varLambda_\alpha$ is
obtained by inducing the character $\alpha$ of the group $A$ to
$G$. In the  proof of Ramsey and Weit, it is crucial that the
function $h$ appearing in the factorization of $f$ commutes with
$\mu$.  In a motion group $G=A\times_\varphi K$ however, the
center of $M(G)$ may not contain non-trivial elements of $L^1(G)$.
Yet, for certain $f\in L^1_0(G)$, we are able to exhibit
appropriate measures $\nu$ in the center of $M(G)$ which may be
used in the place of $h$ in the above argument. Then we use the
spectral radius formula for measures \eqref{SRFM} to conclude that
such $f$ satisfy the mixing condition if $\mu$ is spread-out and
satisfies \eqref{SR}.  By a result of Ludwig on sets of spectral
synthesis we then obtain that such $f$ are dense in
$\ker(\varLambda_0)$. However, since $\ker(\varLambda_0)$ may be
strictly contained in $ L^1_0(G)$, an additional argument is
required in order to treat the general $f\in L^1_0(G)$.

Finally, we mention two more papers that are related. In
\cite{Kan} Kaniuth considers more general groups $G$, namely
locally compact Hausdorff groups of polynomial growth and with a
symmetric group-algebra $L^1(G)$, but only \emph{central\/}
probability measures $\mu\in M(G)$ on such groups; for such
measures he gives the necessary and sufficient conditions
$\norm{\hmu(U)}<1$ and $\hmu(U)\neq I$, for all non-trivial
irreducible $U$, for $\mu$ to be mixing and ergodic by
convolutions respectively. Also related, although more loosely, is
the paper by Jones, Rosenblatt and Tempelman \cite{JRT}, which,
however, has a wider scope.

\bigskip

\noindent We close this Section by fixing some notation and
recalling some terminology, to be used throughout the paper.

\begin{notterm}
We shall follow the terminology of \cite{Fol} regarding
group-representations.  In particular, by a unitary representation
of a locally compact Hausdorff group $G$ we shall always mean a
group homomorphism from $G$ into the group of unitary operators on
some Hilbert space, which is continuous with respect to the strong
operator topology.  Irreducible will always mean topologically
irreducible.  Recall also that the unitary dual $\hG$ of $G$
consists of unitary equivalence classes of irreducible unitary
representations of $G$; for such a representation $U$, we shall
denote by $[U]$ the equivalence class in $\hG$ to which $U$
belongs, by $\mathcal{H}_U$ the representation space of $U$, and
by $d_{[U]}$ the dimension of $\mathcal{H}_U$.

Let $G$ be a locally compact Hausdorff group.  We shall denote by
$M(G)$ the Banach-$\ast$ algebra of complex, regular Borel
measures on $G$.  $L^1(G)$ will stand for the sub-algebra of
$M(G)$ consisting of Haar-integrable Borel functions on $G$ and
$L^1_0(G)$ for the closed two-sided ideal of $f\in L^1(G)$ with
$\int_G f\,d\lambda_G=0$.  For a measure $\mu\in M(G)$, the
Fourier transform of $\mu$ is the bounded linear operator
$\hmu(U):=\int_G U(x^{-1})\, d\mu(x)$, defined, for any continuous
unitary representation $U$ of $G$ on some Hilbert space
$\mathcal{H}_U$, weakly by
$\langle\hmu(U)u,v\rangle:=\int_G\langle
U(x^{-1})u,v\rangle\,d\mu(x)$ $(u,v\in\mathcal{H}_U)$. For such a
representation $U$, we shall also write $U(\mu):=\int_G U(x)\,
d\mu(x)$ for the $*$-representation that $U$ induces on $M(G)$;
notice that $\hmu(U)=U(\bar{\mu})^*$, where $\mbox{}^*$ denotes
adjoint and $\bar{\mu}$ complex conjugation:
$\bar{\mu}(B)=\overline{\mu(B)}$ for Borel subsets $B$ of $G$.

If $T$ is a bounded linear operator on a Hilbert or Banach space,
we shall denote by $\sigma(T)$ its spectrum and by $\varrho(T)$
its spectral radius.

If $\mathcal{H}$ is a Hilbert space, $\boldsymbol{B}(\mathcal{H})$
shall denote the space of bounded linear operators from
$\mathcal{H}$ to $\mathcal{H}$.

If $E$ is a set in a space $X$, we shall denote by $\mathbf{1}_E$
the function which is $1$ on $E$ and $0$ elsewhere; thus, in
particular, if $G$ is a locally compact Hausdorff group,
$\mathbf{1}_G$ identifies with the trivial representation of $G$.

When $G$ is compact, we will always assume that $\lambda_G$ has
been chosen to satisfy $\lambda_G(G)=1$.

Finally, all groups considered in the paper will be assumed to
have Hausdorff topologies, without further notice.
\end{notterm}

\section{A Necessary Condition for Mixing in General Groups}
\label{nec-cond}

In this Section we give some necessary conditions for mixing by
convolutions for general groups $G$. We shall use the notion of a
quasi-compact operator, and recall the definition here.  This
class of operators was introduced by Kryloff and Bogolio\`{u}boff
(\cite{KB1,KB2}).

\begin{definition} A linear operator $T$ on a Banach space $X$ is
\emph{quasi-compact} if there exist $n\in\mathbb{N}$ and a compact
operator $Q$ on $X$ such that $\norm{T^n-Q}<1$.
\end{definition}

We single out the following property of quasi-compact operators
(see the Remarks following Theorem 2.2.8 and the discussion
following Theorem 2.2.7 of \cite{Kr}).

\begin{lemma}[Yosida--Kakutani \cite{YoKa}]
\label{Yosida-Kakutani-1} Let $T$ be a quasi-compact operator on a
Banach space $X$ such that
$\sup_{n\in\mathbb{N}}\norm{T^n}<\infty$.  Then either
$\varrho(T)<1$ or $\set{z\in\sigma(T)\colon\abs{z}=1}$ contains
only eigenvalues of $T$.
\end{lemma}

The following Lemma sheds then some light on the role played by
the spread-out condition:

\begin{lemma}
\label{lem1} Let $G$ be a locally compact \textup{CCR} group. If
$\mu$ is a spread-out probability measure in $M(G)$\textup{,} then
$\hmu(U)$ is quasi-compact for any $[U]\in\hG$.
\end{lemma}

\begin{proof}
Write $\mu^n=(\mu^n)_{\mathrm{a.c.}}+(\mu^n)_{\mathrm{s}}$ with
$(\mu^n)_{\mathrm{s}}\perp\lambda_G$ and
$(\mu^n)_{\mathrm{a.c.}}\ll\lambda_G$, $n\in\mathbb{N}$.  If
$\mu^n\neq (\mu^n)_{\mathrm{s}}$ for some $n$, then
\[\Bignorm{\hmu(U)^n-\widehat{(\mu^n)_{\mathrm{a.c.}}}(U)} =
\Bignorm{\widehat{(\mu^n)_{\mathrm{s}}}(U)} \ls
\norm{(\mu^n)_{\mathrm{s}}}<1 ,\] and
$\widehat{(\mu^n)_{\mathrm{a.c.}}}(U)$ is compact.
\end{proof}

We shall also use the following fact:

\begin{lemma}\label{lem2}
Let $G$ be a locally compact group.  Then\textup{,} for any
$[U]\in\hG\smallsetminus\set{\mathbf{1}_G}$\textup{,} there exists
$h\in\mathcal{H}_U$ such that $\set{U(f)h\colon f\in L_0^1(G)}$ is
dense in $\mathcal{H}_U$.
\end{lemma}

\begin{proof}  For any $h\neq 0$, the subspace
$\set{U(f)h\colon f\in L_0^1(G)}$ of $\mathcal{H}_U$ is invariant
under $U$, and hence it is either $\set{0}$ or dense; thus, we
only have to exclude the possibility that it is $\set{0}$ for
every $h$.  Since $L^1_0(G)$ has co-dimension one in $L^1(G)$, if
these subspaces are all trivial, then $U$ is one-dimensional, and
the representation $f\mapsto U(f)$ of $L^1(G)$ has kernel
$L^1_0(G)$. But the only one-dimensional representation of $G$ for
which the corresponding representation of $L^1(G)$ has kernel
$L^1_0(G)$ is $\mathbf{1}_G$.
\end{proof}

\begin{proposition}
\label{prop2} Let $G$ be a locally compact group\textup{,} and let
$\mu$ be a probability measure in $M(G)$ which is mixing by
convolutions. Then\textup{:}
\begin{enumerate}
\item[\1]
$\hmu(U)^n\to 0$ in the strong operator topology\textup{,} for any
$[U]\in\hG\smallsetminus\set{\mathbf{1}_G}$.
\item[\2]
$\varrho\bigl(\hmu(U)\bigr)<1$ for any
$[U]\in\hG\smallsetminus\set{\mathbf{1}_G}$ for which $\hmu(U)$ is
quasi-compact.
\end{enumerate}
\end{proposition}

\begin{proof} \1\
Let $[U]\in\hG\smallsetminus\set{\mathbf{1}_G}$ and fix
$h\in\mathcal{H}_U$. Let also $\epsilon>0$ be given.  By Lemma
\ref{lem2} there exists $f\in L_0^1(U)$ and $h'\in\mathcal{H}_U$
such that $\norm{U(f)h'-h}<\epsilon/2$.  Set
$g(x):=\varDelta_G(x^{-1}) f(x^{-1})$, where $\varDelta_G$ is the
modular function of $G$. Then $\hg(U)=U(f)$, whence
\[ \norm{\hmu(U)^n h} \ls
\norm{\hmu(U)^n \hg(U) h'} + \norm{\hmu(U)^n (\hg(U) h' - h)} <
\norm{g\ast\mu^n}_1 + \epsilon/2 ,\] and this is $<\epsilon$ for
all sufficiently large $n$, since $g\in L^1_0(G)$.

\smallskip

\noindent\2\ Fix $[U]\in\hG\smallsetminus\set{\mathbf{1}_G}$ and
suppose that $\hmu(U)$ is quasi-compact.  By \1, $\hmu(U)$ cannot
have eigenvalues of modulus one.  It then follows from Lemma
\ref{Yosida-Kakutani-1} that $\varrho\bigl(\hmu(U)\bigr)<1$.
\end{proof}

Note that as a result of Proposition \ref{prop2} and Lemma
\ref{lem1}, one immediately obtains the necessity of condition
\eqref{SR} for mixing of spread-out measures on CCR groups:

\begin{corollary}
\label{cor1} Let $G$ be a locally compact \textup{CCR} group. If
$\mu$ is a spread-out probability measure in $M(G)$ which is
mixing by convolutions\textup{,} then
$\varrho\bigl(\hmu(U)\bigr)<1$ for all
$[U]\in\hG\smallsetminus\set{\mathbf{1}_G}$.
\end{corollary}

Corollary \ref{cor1} will in fact be subsumed by the stronger
result of Proposition \ref{wmix}.

\section{Unitary Representations of Motion Groups}
\label{prelims}

\subsection{Locally Compact Motion Groups}
Let $K$ be a compact group, $A$ an Abelian group, and consider the
semi-direct product $G=A\times_\varphi K$.  So $\varphi\colon K\to
\operatorname{Aut}(A)$ is assumed to be a group homomorphism, and
we shall write $\varphi_\kappa\in\operatorname{Aut}(A)$ for the
image of the element $\kappa\in K$ under $\varphi$. The group
operation on $G$ is given by
\[ (a_1,\kappa_1)\cdot (a_2,\kappa_2) :=
\bigl(a_1+\varphi_{\kappa_1}(a_2),\kappa_1\kappa_2\bigr),\] and
the mapping $(a,\kappa)\mapsto\varphi_\kappa(a)$ is assumed to be
continuous.

Left Haar measure $\lambda_G$ on $G$ is the product
$\lambda_A\otimes\lambda_K$ of the Haar measures on $A$ and $K$
respectively, owing to the fact that we assume $K$ compact and $A$
abelian \cite[15.29]{HeRo1}. Furthermore, using standard
arguments, it is not hard to see that each $\varphi_\kappa$ must
be measure-preserving; i.e, if
$\varphi_\kappa(\lambda_A)=\lambda_A\circ\varphi_\kappa^{-1}$
denotes the measure defined by
$\varphi_\kappa(\lambda_A)(B):=\lambda_A(\varphi_\kappa^{-1}(B))$
for Borel subsets $B$ of $A$, then \begin{equation}\label{inv}
\varphi_\kappa(\lambda_A)=\lambda_A \qquad\forall\, \kappa\in K
\end{equation}
(see \cite[15.29]{HeRo1} again).

Throughout, we shall be using additive notation for Abelian groups
and multiplicative notation for other groups. $1$ shall denote the
neutral element of the group $K$ and $0$ the neutral element of
$A$.  Finally, if $a\in A$ and $\alpha\in\hA$ is a character of
the Abelian group $A$, we shall also use the notation
\[\langle a,\alpha\rangle := \alpha(a). \]

\subsection{Unitary Representations of Motion Groups}
\label{UR} The action $\varphi$ of the group $K$ on $A$ determines
an action of $K$ on the dual group $\hA$ of $A$ through
$\alpha\mapsto
\varphi_\kappa(\alpha):=\alpha\circ\varphi_{\kappa^{-1}}$.

For $\alpha\in\hA$, we shall denote by $\varLambda_\alpha$ the
induced representation $\operatorname{ind}^G_{A\times\set{1}}
(\alpha)$ on $G$.  This may be realized on $L^2(K)$ as follows.
For $(a,\kappa)\in G$ and $\phi\in L^2(K)$,
\begin{equation}\label{varLambda}[\varLambda_\alpha (a,\kappa) \phi](\kappa') = \langle a,
\varphi_{\kappa'}(\alpha)\rangle\cdot [L_K(\kappa) \phi] (\kappa')
= \langle a, \varphi_{\kappa'}(\alpha)\rangle\cdot
\phi(\kappa^{-1}\kappa'),\end{equation} $\kappa'\in K$, where
$L_K$ denotes the left regular representation of $K$ on $L^2(K)$.
Observe that $\varLambda_\alpha$ and $\varLambda_{\alpha'}$ are
unitarily equivalent when $\alpha$ and $\alpha'$ belong to the
same orbit, i.e., when $\alpha'=\varphi_{\kappa'}(\alpha)$ for
some $\kappa'\in K$; in fact
\[
\varLambda_{\varphi_{\kappa'}(\alpha)}(a,\kappa) = R_K(\kappa')\,
\varLambda_\alpha(a,\kappa)\, R_K(\kappa')^{-1}\qquad
((a,\kappa)\in G) \] for all $\alpha\in\hA,\ \kappa^\prime\in K$,
where $R_K$ denotes the right regular representation of $K$ on
$L^2(K)$.

In all our results concerning motion groups, we shall be assuming
that $G$ acts regularly on the dual group $\hA$ of $A$. The
irreducible unitary representations of $G=A\times_\varphi K$ are
then as follows.

\begin{proposition}
\label{prop-1} Let $G=A\times_{\varphi} K$ be a motion group with
$G$ acting regularly on $\hA$.  Then any irreducible unitary
representation of $G$ is unitarily equivalent to a
sub-representation of $\varLambda_\alpha$ for some $\alpha\in\hA$.
Furthermore\textup{,} each $\varLambda_\alpha$ is the direct sum
of irreducible unitary representations of $G$. Finally\textup{,}
two irreducible unitary representations $U_1,U_2$ of $G$ are
unitarily equivalent only if they are equivalent to
sub-representations of $\varLambda_{\alpha_1}$ and
$\varLambda_{\alpha_2}$ respectively\textup{,} with $\alpha_1$ and
$\alpha_2$ belonging to the same orbit\textup{,} i.e.\textup{,}
with $\alpha_2=\varphi_\kappa(\alpha_1)$ for some $\kappa\in K$.
\end{proposition}

\begin{proof}
This follows from \cite[Theorem 6.42]{Fol} and induction by stages
\cite[Theorem 6.14]{Fol}.
\end{proof}

We shall also need an analogue of the Riemann--Lebesgue lemma; as
we were unable to locate the following version in the literature,
we give a proof in the Appendix.

\begin{theorem}[Riemann--Lebesgue Lemma]
\label{app-thm-1} Let $G=A\times_{\varphi} K$ be a motion group
with $G$ acting regularly on $\hA$.  Then\textup{,} for any $f\in
L^1(G)$\textup{,} one has the
following\textup{:}\begin{itemize}\item[\1] Given
$\epsilon>0$\textup{,} there exists a compact set
$\hC\subseteq\hA$\textup{,} such that
$\bigl\lVert\hf(\varLambda_\alpha)\bigr\rVert<\epsilon$ for all
$\alpha\in\hA\smallsetminus\hC$.
\item[\2]
Given $\epsilon>0$ and $\alpha\in\hA$\textup{,} if
$\varLambda_\alpha = \bigoplus_{i\in\mathcal{I}} U_i$ is a
direct-sum decomposition of $\varLambda_\alpha$ into irreducible
unitary representations of $G$, then
$\bigl\lVert\hf(U_i)\bigr\rVert < \epsilon$ for all but finitely
many $i\in\mathcal{I}$.
\end{itemize}
\end{theorem}

Finally, we shall also use the following result, a proof of which
is also given in the Appendix.

\begin{theorem}
\label{app-thm-2} Let $G=A\times_{\varphi} K$ be a motion group
with $G$ acting regularly on $\hA$.  Then\textup{,} for any
$\mu\in M(G)$\textup{,} the operator-valued function
$\alpha\mapsto\hmu(\varLambda_\alpha)$ is uniformly continuous on
$\hA$ with respect to the norm topology on
$\boldsymbol{B}(L^2(K))$.
\end{theorem}

\section{Spectral Radius Formulae in Motion Groups}
\label{SpecRad}

\subsection{The Analogue of the Group $C^*$-Algebra for
Measures}\label{C^*}

We shall need to consider the analogue of the group-$C^*$-algebra
$C^*(G)$ for measures on $G$. Let $G$ be an arbitrary locally
compact group.  For $\mu\in M(G)$, define
\[\norm{\mu}_*:=\sup_{[U]\in\hG}\norm{U(\mu)},\]
where  $U(\mu):=\int_G U(x)\, d\mu(x)$. One verifies easily that
$\mu\mapsto\norm{\mu}_\ast$ is a norm on $M(G)$ (the implication
$\norm{\mu}_\ast=0\,\Longrightarrow\, \mu=0$ follows from the
injectivity of the Fourier transform \cite[18.2.3]{Dix}). We shall
then denote by $D^*(G)$ the completion of the unital Banach
algebra $M(G)$ with respect to this norm. Then, $D^*(G)$ is a
unital $C^\ast$-algebra, and the group-$C^\ast$-algebra
$C^\ast(G)$ is a closed sub-algebra of $D^*(G)$. For compact
groups, $D^*(G)$ has also been considered in \cite{AG}.

Any $*$-representation of $M(G)$ extends uniquely to a
$*$-representation of $D^*(G)$, so in particular, if $U$ is any
irreducible unitary representation of $G$, $\mu\mapsto U(\mu)$
extends to a $*$-representation of $D^*(G)$. The Fourier transform
$\mu\mapsto\hmu(U)$ then also extends to $D^*(G)$.

Finally, we shall also use the following fact (see \cite[Theorem
22.11]{HeRo1}).

\begin{fact}\label{GF}
Let $G$ be a locally compact group and denote by $L_G$ its left
regular representation on $L^2(G)$.  Then
$\hmu(L_G)=0\,\Longrightarrow\, \mu=0$.
\end{fact}

\subsection{Spectral Radius Formulae}

Recall that the spectral radius of an element $a$ in a Banach
algebra $\mathcal{A}$ may be defined by
\begin{equation}\label{sr} \varrho(a) =
\lim_{n\to\infty}\norm{a^n}^{1/n} = \inf\limits_{n\in\mathbb{N}}
\norm{a^n}^{1/n} .\end{equation}  One then has the following
spectral radius formula for measures in motion groups:

\begin{theorem}
\label{srfm} Let $G=A\times_{\varphi} K$ be a motion group with
$G$ acting regularly on $\hA$. Then\textup{,} for any $\mu\in
M(G)$\textup{,} one has that
\begin{align*}
\varrho(\mu)=\lim\limits_{n\to\infty}\norm{\mu^n}^{1/n} &=
\sup\limits_{[U]\in\hG}\varrho\bigl(\hmu(U)\bigr) \vee
\inf\limits_{n\in\mathbb{N}} \norm{(\mu^n)_{\mathrm{s}}}^{1/n} \\
&=
\sup\limits_{\alpha\in\hA}\varrho\bigl(\hmu(\varLambda_\alpha)\bigr)
\vee \inf\limits_{n\in\mathbb{N}}
\norm{(\mu^n)_{\mathrm{s}}}^{1/n} ,\end{align*} where
$\varrho(\mu)$ denotes the spectral radius of $\mu$ in the unital
Banach algebra $M(G)$\textup{,} and $(\mu^n)_{\mathrm{s}}$ the
singular part of $\mu^n$ with respect to Haar measure $\lambda_G$.
\end{theorem}

\begin{note*} For numbers $a,b\in\mathbb{R}$, $a\vee
b:=\max\set{a,b}$.\end{note*}

\begin{lemma}
\label{srfl} Let $G=A\times_{\varphi} K$ be a motion group with
$G$ acting regularly on $\hA$.  Then\textup{,} for any $\mu\in
M(G)$\textup{,} \[ \lim\limits_{n\to\infty} \norm{\mu^n}^{1/n} \ls
\norm{\mu_\mathrm{s}}\vee \sup\limits_{[U]\in\hG}
\varrho\bigl(\hmu(U)\bigr) = \norm{\mu_\mathrm{s}} \vee
\sup\limits_{\alpha\in\hA}\varrho\bigl(\hmu(\varLambda_\alpha)\bigr)
.\]
\end{lemma}

\begin{proof}[Proof of Theorem \textup{\ref{srfm}}]
Let $\mu\in M(G)$.  Since 
\[\norm{\mu^n}\gr\norm{\hmu(U)^n}\gr \left[
\varrho\bigl(\hmu(U)\bigr) \right]^n \] for all $[U]\in\hG$ and
$n\in\mathbb{N}$, it is clear that $\varrho(\mu)\gr
\sup_{[U]\in\hG} \varrho\bigl(\hmu(U)\bigr)$. Since also \[
\norm{\mu^n} = \norm{(\mu^n)_{\mathrm{a.c.}}} +
\norm{(\mu^n)_{\mathrm{s}}} \gr \norm{(\mu^n)_{\mathrm{s}}} \]
(cf.\ \cite[Theorem 14.22]{HeRo1} for the equality), it follows
that
\begin{equation}\label{srfm1} \varrho(\mu) \gr
\sup\limits_{[U]\in\hG} \varrho\bigl(\hmu(U)\bigr) \vee
\inf\limits_{n\in\mathbb{N}} \norm{(\mu^n)_{\mathrm{s}}}^{1/n} .
\end{equation}

For the reverse inequality apply Lemma \ref{srfl} to the powers
$\mu^n$ of $\mu$.  For each $n\in\mathbb{N}$, \[ \varrho(\mu^n)\ls
\norm{(\mu^n)_{\mathrm{s}}}\vee \sup\limits_{[U]\in\hG}
\varrho\bigl( \hmu(U)^n\bigr) ,\] and since in any Banach algebra
$\varrho(a^m)=[\varrho(a)]^m$ for all $a$ and $m$, it follows that
\[ [\varrho(\mu)]^n = \varrho(\mu^n) \ls
\norm{(\mu^n)_{\mathrm{s}}}\vee \sup\limits_{[U]\in\hG}
\varrho\bigl( \hmu(U)^n\bigr) .\]  This, together with
\eqref{srfm1} show that \[\label{srfm3} \varrho(\mu) =
\sup\limits_{[U]\in\hG} \varrho\bigl(\hmu(U)\bigr) \vee
\inf\limits_{n\in\mathbb{N}} \norm{(\mu^n)_{\mathrm{s}}}^{1/n} .\]

The equality
\[\varrho(\mu) =
\sup\limits_{\alpha\in\hA}\varrho\bigl(\hmu(\varLambda_\alpha)\bigr)
\vee \inf\limits_{n\in\mathbb{N}}
\norm{(\mu^n)_{\mathrm{s}}}^{1/n} \] is established in the same
way.
\end{proof}

\begin{proof}[Proof of Lemma \textup{\ref{srfl}}]
As $L^1(G)$ is a symmetric Banach $^*$-algebra \cite{Gang}, Remark
3 of \cite{PaA} yields the formula
\begin{equation}\label{srfl0}\varrho(f) =
\lim\limits_{n\to\infty}\norm{f^n}_1^{1/n} =
\lim\limits_{n\to\infty}\norm{f^n}_*^{1/n}
\end{equation} for $f\in L^1(G)$, where $\norm{\ }_*$ is the norm
\[ \norm{\mu}_* := \sup_{[U]\in\hG} \norm{\hmu(U)} \] on $M(G)$,
defined in Subsection \ref{C^*}. Write also \[ \varrho_*(\mu) :=
\lim\limits_{n\to\infty}\norm{\mu^n}_*^{1/n} \qquad (\mu\in
M(G)),\] and recall that $D^*(G)$ is the completion of $M(G)$ with
respect to the norm $\norm{\ }_*$.  Since $D^*(G)$ and $M(G)$ are
unital Banach algebras, we have that \[ \varrho_*(\mu) =
\sup\set{\abs{z}\colon z\in\mathbb{C} \mbox{ and } z\delta_{e} -
\mu \mbox{ is not invertible in } D^*(G)} \] and
\[ \varrho(\mu) = \sup\set{\abs{z}\colon
z\in\mathbb{C} \mbox{ and } z\delta_{e} - \mu \mbox{ is not
invertible in } M(G)} \] for $\mu\in M(G)$, where $e$ is the
neutral element in $G$ and $\delta_e$ the Dirac point-mass at $e$.
We shall first show that
\begin{equation}\label{srfl1}
\varrho(\mu)\ls\varrho_*(\mu)\vee\varrho(\mu_{\mathrm{s}})
\end{equation} for any $\mu\in M(G)$.

Indeed, fix $\mu\in M(G)$.  If $z\in\mathbb{C}$ with $\abs{z} >
\varrho_*(\mu)\vee\varrho(\mu_{\mathrm{s}})$, then
$\delta_e-z^{-1}\mu$ is invertible in $D^*(G)$ and
$\delta_e-z^{-1}\mu_\mathrm{s}$ is invertible in $M(G)$. Set
\[\nu:=\left(\delta_e-
z^{-1}\mu_\mathrm{s}\right)^{-1}*\mu_\mathrm{a.c.} ,\] and notice
that, because $L^1(G)$ is an ideal in $M(G)$, $\nu\in L^1(G)$.
Since
\begin{equation}\label{srfl2}
\left(\delta_e-z^{-1}\mu_\mathrm{s}\right)^{-1}*
\left(\delta_e-z^{-1}\mu\right)= \delta_e-z^{-1}\nu ,
\end{equation}
the right-hand side is invertible in $D^*(G)$, because the left
side is.  This shows that $\varrho_*(\nu)\ls
\varrho_*(\mu)\vee\varrho(\mu_{\mathrm{s}})$.  By \eqref{srfl0}
then, $\varrho(\nu)=\varrho_*(\nu)\ls
\varrho_*(\mu)\vee\varrho(\mu_{\mathrm{s}})$, since $\nu\in
L^1(G)$.  Thus if $z\in\mathbb{C}$ with $\abs{z} >
\varrho_*(\mu)\vee\varrho(\mu_{\mathrm{s}})$, then $\abs{z} >
\varrho(\nu)$ and the inverse $\left(\delta_e-z^{-1}\nu
\right)^{-1}$, whose existence in $D^*(G)$ is guaranteed by
\eqref{srfl2}, must actually belong to $M(G)$.  Now \eqref{srfl2}
also yields that \[ \left(\delta_e-z^{-1}\mu\right)^{-1}=
\left(\delta_e-z^{-1}\nu\right)^{-1} *
\left(\delta_e-z^{-1}\mu_s\right)^{-1} ,\] whence
$\left(\delta_e-z^{-1}\mu\right)^{-1}\in M(G) * M(G) = M(G)$. This
shows \eqref{srfl1}.

In order to establish the Lemma, we now only need to show that
\begin{equation}
\label{srfl3} \varrho_*(\mu) \ls
\sup_{[U]\in\hG}\varrho\bigl(\hmu(U)\bigr) \vee
\norm{\mu_{\mathrm{s}}} =
\sup_{\alpha\in\hA}\varrho\bigl(\hmu(\varLambda_\alpha)\bigr) \vee
\norm{\mu_{\mathrm{s}}}.
\end{equation}
Let $\mu_{\mathrm{a.c.}}$ denote the absolutely continuous part of
$\mu$, and let $\epsilon >0$ be given. By the Riemann--Lebesgue
lemma (Theorem \ref{app-thm-1}), there exists a compact set
$\hC\subseteq\hA$ such that
$\norm{\hmu_{\mathrm{a.c.}}(\varLambda_\alpha)}<\epsilon$ for all
$\alpha$ in $\hA\smallsetminus\hC$.

Next set
\[ r_n(\alpha):=
\norm{\hmu(\varLambda_\alpha)^n}^{1/n} \qquad\mbox{and}\qquad
r(\alpha):= \varrho\bigl(\hmu(\varLambda_\alpha)\bigr) .\] The
inequality
\begin{align*} \norm{\hmu(\varLambda_\alpha)^n -
\hmu(\varLambda_\beta)^n} &\ls \norm{\hmu(\varLambda_\alpha) -
\hmu(\varLambda_\beta)}\, \sum\limits_{k=0}^{n-1}
\norm{\hmu(\varLambda_\alpha)^k}\,
\norm{\hmu(\varLambda_\beta)}^{n-k-1}
\\ &\ls n\, \norm{\mu}^{n-1}\, \norm{\hmu(\varLambda_\alpha) -
\hmu(\varLambda_\beta)} ,
\end{align*}
together with the norm-continuity of the operator-valued function
$\alpha\mapsto\hmu(\varLambda_\alpha)$ (Theorem \ref{app-thm-2}),
show that each $r_n$ is a continuous function on $\hA$. The
norm-continuity of $\alpha\mapsto\hmu(\varLambda_\alpha)$ also
implies that $r$ is upper semi-continuous; hence it attains its
maximum on $\hC$, and let $r^*=\max_{\alpha\in\hC} r(\alpha)$.
Since the mapping $\alpha\mapsto r_n(\alpha)-r^*$ is continuous,
the sets $\hC_n:=\bigset{\alpha\in\hC\colon
r_n(\alpha)-r^*\gr\epsilon}$ are compact, and since
$r_n(\alpha)\downarrow r(\alpha)\ls r^*$ for each $\alpha$, by
\eqref{sr}, $\bigcap_{n=1}^\infty \hC_n = \varnothing$; it follows
that for some $n(\epsilon)\in\mathbb{N}$,
$\hC_{n(\epsilon)}=\varnothing$.  Then, for $n\gr n(\epsilon)$,
one has that
\begin{align}\label{SR-3}
\sup\limits_{\alpha\in\hA}\norm{\hmu(\varLambda_\alpha)^n} &\ls
\sup\limits_{\alpha\in\hC}\norm{\hmu(\varLambda_\alpha)^n}\vee
\sup\limits_{\alpha\in\hA\smallsetminus\hC}\norm{\hmu(\varLambda_\alpha)}^n\\
&\ls
\sup\limits_{\alpha\in\hC}\norm{\hmu(\varLambda_\alpha)^n}\vee
\sup\limits_{\alpha\in\hA\smallsetminus\hC}
(\norm{\widehat{\mu_{\mathrm{a.c.}}}(\varLambda_\alpha)}+
\norm{\widehat{\mu_{\mathrm{s}}}(\varLambda_\alpha)})^n\nonumber\\
&\ls
\sup\limits_{\alpha\in\hA}[\varrho\bigl(\hmu(\varLambda_\alpha)\bigr)+\epsilon]^n
\vee (\epsilon + \norm{\mu_s})^n\nonumber .
\end{align}
Since, by Proposition \ref{prop-1}, any irreducible unitary
representation $U$ of $G$ is unitarily equivalent to a
sub-representation of $\varLambda_\alpha$ for some $\alpha\in\hA$,
\eqref{SR-3} then implies that
\[ \sup\limits_{[U]\in\hG} \norm{\hmu(U)^n} \ls
\sup\limits_{\alpha\in\hA}\norm{\hmu(\varLambda_\alpha)^n} \ls
\sup\limits_{\alpha\in\hA}[\varrho\bigl(\hmu(\varLambda_\alpha)\bigr)+\epsilon]^n
\vee (\epsilon + \norm{\mu_s})^n, \] whence
\[ \norm{\mu^n}_\ast\ls \sup\limits_{\alpha\in\hA}[\varrho\bigl(\hmu(\varLambda_\alpha)\bigr)+\epsilon]^n
\vee (\epsilon + \norm{\mu_s})^n,\] for all $n\gr n(\epsilon)$.
Since $\epsilon $ was arbitrary, this shows that
\[ \varrho_*(\mu)\ls \sup\limits_{\alpha\in\hA} \varrho\bigl(\hmu(\varLambda_\alpha)\bigr) \vee
\norm{\mu_{\mathrm{s}}} .\] We shall next show that
\begin{equation}\label{SR-4} \sup\limits_{\alpha\in\hA}
\varrho\bigl(\hmu(\varLambda_\alpha)\bigr) \ls
\sup\limits_{[U]\in\hG} \varrho\bigl(\hmu(U)\bigr) \vee
\norm{\mu_{\mathrm{s}}} ,\end{equation} which will show
\eqref{srfl3} and thus complete the proof.

Fix $\alpha\in\hA$.  By Proposition \ref{prop-1},
$\varLambda_\alpha$ is a direct sum of irreducible unitary
representations of $G$, say $\varLambda_\alpha =
\bigoplus_{i\in\mathcal{I}} U_i$.  Then, for any $n\in\mathbb{N}$,
\begin{equation}\label{SR-5} \left[\varrho\bigl(\hmu(\varLambda_\alpha)\bigr)\right]^n \ls
\norm{\hmu(\varLambda_\alpha)^n} = \sup\limits_{i\in\mathcal{I}}
\norm{\hmu(U_i)^n} .\end{equation}  Let $\epsilon >0$.  By Theorem
\ref{app-thm-1}, there exists a finite set
$\mathcal{I}_1\subseteq\mathcal{I}$ such that
$\norm{\widehat{\mu_{\mathrm{a.c.}}}(U_i)} < \epsilon$ for all
$i\in\mathcal{I}_2 := \mathcal{I}\smallsetminus\mathcal{I}_1$.
Choose $n(\epsilon)\in\mathbb{N}$ so that $\norm{\hmu(U_i)^n} \ls
\left[\varrho\bigl(\hmu(U_i)\bigr) + \epsilon\right]^n$ for all
$n\gr n(\epsilon)$ and all $i\in\mathcal{I}_1$.  Then, for $n\gr
n(\epsilon)$, one has that
\begin{align}\label{SR-6}
\sup\limits_{i\in\mathcal{I}} \norm{\hmu(U_i)^n} &\ls
\sup\limits_{i\in\mathcal{I}_1} \norm{\hmu(U_i)^n} \vee
\sup\limits_{i\in\mathcal{I}_2} \norm{\hmu(U_i)^n} \\ &\ls
\sup\limits_{i\in\mathcal{I}_1} \left[\varrho\bigl(\hmu(U_i)\bigr)
+ \epsilon\right]^n \vee \sup\limits_{i\in\mathcal{I}_2}
(\norm{\widehat{\mu_{\mathrm{a.c.}}}(U_i)} +
\norm{\widehat{\mu_{\mathrm{s}}}(U_i)})^n\nonumber
\\ &\ls \sup\limits_{[U]\in\hG} \left[\varrho\bigl(\hmu(U)\bigr) +
\epsilon\right]^n \vee (\epsilon +
\norm{\mu_{\mathrm{s}}})^n,\nonumber
\end{align}
whence by \eqref{SR-5}
\[\varrho\bigl(\hmu(\varLambda_\alpha)\bigr)
\ls \sup\limits_{[U]\in\hG} \left[\varrho\bigl(\hmu(U)\bigr) +
\epsilon\right] \vee (\epsilon + \norm{\mu_{\mathrm{s}}}) .\] As
$\alpha$ and $\epsilon$ were arbitrary, this establishes
\eqref{SR-4}.
\end{proof}

We close this Section with a result from \cite{Pal2} that we shall
use in the sequel (Section \ref{ergodic}):

\begin{proposition}\label{spectra}
Suppose that $G=A\times_{\varphi} K$ is a motion group with $G$
acting regularly on $\hA$.  For $f\in L^1(G),$ let
$\sigma_{M(G)}(f)$ denote the spectrum of $f$ as an element of the
unital Banach algebra $M(G)$ and $\sigma_{D^\ast(G)}(f)$ the
spectrum of $f$ as an element of the unital Banach algebra
$D^\ast(G)$.  Then\textup{:}
\[\sigma_{D^\ast(G)}(f)\cup\set{0}=\sigma_{M(G)}(f)\cup\set{0}\quad
\forall\, f\in L^1(G).\]
\end{proposition}

\begin{proof}  This follows from \cite[Proposition
10.4.6\,(b)]{Pal2}.
\end{proof}

\section{Mixing in Motion Groups}

In this Section we prove the following result:

\begin{theorem}\label{main}
Let $G=A\times_{\varphi} K$ be a motion group with $G$ acting
regularly on $\hA$. Then\textup{,} if $\mu\in M(G)$ is a
spread-out probability measure\textup{,} $\mu$ is mixing by
convolutions if and only if $\varrho\bigl(\hmu(U)\bigr)< 1\
\forall\, [U]\in\hG\smallsetminus\set{\mathbf{1}_G}$.
\end{theorem}

Observe that one direction of the Theorem follows directly from
Corollary \ref{cor1}, as any motion group with $G$ acting
regularly on $\hA$ is CCR  \cite[4.5.2.1]{War}.

For the converse, first recall the standard facts that any motion
group has a) polynomial growth \cite[Theorem 1.4]{Guiv}, and b) a
symmetric group algebra $L^1(G)$ \cite{Gang}.

Given a probability measure $\mu\in M(G)$, set
\begin{equation}\label{Im} I_\mu:=\set{f\in L^1(G)\colon
\norm{f\ast\mu^n}_1\to 0\ \text{as}\ n\to\infty};\end{equation}
$I_\mu$ is clearly a closed left ideal in $L^1(G)$, contained in
$L^1_0(G)$. Notice however, that it is not \emph{a priori\/} clear
that $I_\mu$ is a two-sided ideal in $L^1(G)$; for this reason,
one cannot directly refer to spectral synthesis, as is the case
when $G$ is an Abelian group (cf.\ \cite{Ram}) or when $\mu$ is a
central measure (cf.\ \cite{Kan}), and we shall have to use Lemma
\ref{lem6-1} below instead.

To prove Theorem \ref{main} we will first show that, if $\mu$ is
spread-out and satisfies $\varrho\bigl(\hmu(U)\bigr)< 1\ \forall\,
[U]\in\hG\smallsetminus\set{\mathbf{1}_G}$, then
$\ker(\varLambda_0)\subseteq I_\mu$, and then deduce from this
that all of $L^1_0(G)$ is contained in $I_\mu$.

Let $\widehat{\mathcal{C}}_0$ denote the collection of all compact
subsets of $\hA$ not containing $0\in\hA$\textup{,} and set
\begin{equation}\label{I}
I:=\bigcup\limits_{\hC\in\widehat{\mathcal{C}}_0}\bigset{f\in
L^1(G)\colon \hf(\varLambda_\alpha)=0\
\forall\,\alpha\in\hA\smallsetminus \hC}.\end{equation}

\begin{lemma}\label{central}
Given a compact subset $\hC$ of $\hA$ not containing $0\in\hA,$
there exists $\nu\in Z(M(G))$ such that
$\hnu(\varLambda_\alpha)=\hh(\alpha)I_{L^2(K)},$ where
$I_{L^2(K)}$ is the identity operator on $L^2(K),$ and where
$h\colon A\to\mathbb{C}$ is a \textup{(}necessarily\textup{)}
$K$-invariant $L^1$-function whose Fourier transform
satisfies\textup{:}
\begin{itemize}\item[\1] $\hh(\alpha)=1\ \forall\,\alpha\in\hC;$
\item[\2] $\hh(\alpha)=0$ for all $\alpha$ outside a compact set
not containing $0\in\hA;$
\item[\3] $0\ls\hh\ls 1$.
\end{itemize}
\end{lemma}

\begin{note*} $Z(M(G))$ denotes the center of $M(G)$, i.e.,
those elements of $M(G)$ which commute with every other element of
$M(G)$.
\end{note*}

\begin{proof}[Proof of Lemma \textup{\ref{central}}]
Have a $K$-invariant function $h\in L^1(A)$ whose Fourier
transform $\hh$ satisfies properties \1--\3, and let $\nu\in M(G)$
be the measure on $G$ defined by
$\nu:=(h\,d\lambda_A)\otimes\delta_{\set{1}}$, where $\otimes$
denotes product-measure and $\delta_{\set{1}}$ denotes a
point-mass at the identity $1\in K$.  Then, by \eqref{varLambda},
\begin{align*}
[\hnu(\varLambda_\alpha) \phi] (\kappa) = \int_A \langle -a,
\varphi_{\kappa}(\alpha)\rangle\, \phi(\kappa)\,
h(a)\,d\lambda_A(a) = \phi(\kappa) \,
\hh(\varphi_{\kappa}(\alpha))
\end{align*} for any $\alpha\in\hA$, $\kappa\in K$, and $\phi\in
L^2(K)$ (recall that the representation space of each
$\varLambda_\xi$ is $L^2(K)$).  Since $h$ and $\lambda_A$ are
$K$-invariant, we have that
$\hh(\varphi_{\tkappa}(\alpha))=\hh(\alpha)$, and so we finally
deduce that \begin{equation}\label{hnu} \hnu(\varLambda_\alpha) =
\hh(\alpha)\, I_{L^2(K)}\qquad (\alpha\in\hA)\end{equation} where
$I_{L^2(K)}$ is the identity operator on $L^2(K)$.
\end{proof}

\begin{lemma}\label{lem6-2}
Let $G$ be as in Theorem \textup{\ref{main}} and let $\mu\in M(G)$
be a spread-out probability measure with
$\varrho\bigl(\hmu(U)\bigr)< 1\ \forall\,
[U]\in\hG\smallsetminus\set{\mathbf{1}_G}$. Then $I\subseteq
I_\mu$.
\end{lemma}

\begin{proof}
Fix $f\in I$.  Then $\hf(\varLambda_\alpha)=0$ for all
$\alpha\in\hA\smallsetminus\hC$ for some compact set
$\hC\in\widehat{\mathcal{C}}_0$, which we fix. Fix a central
measure $\nu$ as in Lemma \ref{central} for this $\hC$, and denote
the support of the corresponding function $\hh$ by $\hS$. Then
$\hf(\varLambda_\alpha)=\hf(\varLambda_\alpha)\,\hnu(\varLambda_\alpha)=
\hnu(\varLambda_\alpha)\,\hf(\varLambda_\alpha)$ for all
$\alpha\in\hA$, and since each irreducible unitary representation
of $G$ is contained as a direct summand in some
$\varLambda_\alpha$, it follows that
$\hf(U)\,\hnu(U)=\hnu(U)\,\hf(U)=\hf(U)$ for any irreducible
unitary representation $U$ of $G$.  Hence $f\ast\nu=\nu\ast f=f$,
by the injectivity of the Fourier transform on $M(G)$ ($\norm{\
}_\ast$, defined in Subsection \ref{C^*}, is a norm on $M(G)$).
Since $\nu$ is also central, we then have that
\begin{equation}\label{final1}
\norm{f\ast\mu^n}_1=\norm{f\ast\nu^n\ast\mu^n}_1=
\norm{f\ast(\nu\ast\mu)^n}_1\ls\norm{f}_1\norm{(\nu\ast\mu)^n}.\end{equation}
On the other hand, by the spectral radius formula of Theorem
\ref{srfm},
\begin{align}\label{final2}
\lim\limits_{n\to\infty}\norm{(\nu\ast\mu)^n}^{1/n} &=
\sup\limits_{\alpha\in\hA}
\varrho\bigl(\hmu(\varLambda_\alpha)\hnu(\varLambda_\alpha)\bigr)\vee
\inf\limits_{n\in\mathbb{N}}\norm{(\nu^n\ast
\mu^n)_{\mathrm{s}}}^{1/n}\\ &=  \sup\limits_{\alpha\in\hA}
\Bigl[\hh(\alpha)\varrho\bigl(\hmu(\varLambda_\alpha)\bigr)\Bigr]\vee
\inf\limits_{n\in\mathbb{N}}\norm{(\nu^n\ast
\mu^n)_{\mathrm{s}}}^{1/n}\nonumber\\ &=
\sup\limits_{\alpha\in\hS}
\Bigl[\hh(\alpha)\varrho\bigl(\hmu(\varLambda_\alpha)\bigr)\Bigr]\vee
\inf\limits_{n\in\mathbb{N}}\norm{(\nu^n\ast
\mu^n)_{\mathrm{s}}}^{1/n},\nonumber\end{align} and this is $<1$
because we are assuming that $\varrho\bigl(\hmu(U)\bigr)<1$ for
all $[U]\in\hG\smallsetminus\set{\mathbf{1}_G},$ and because $\mu$
is spread-out and $0\ls\hh\ls 1$.  Indeed, since
\[\norm{(\nu^n\ast\mu^n)_{\mathrm{s}}}\ls
\norm{(\nu^n)_{\mathrm{s}}\ast(\mu^n)_{\mathrm{s}}}\ls
\norm{(\mu^n)_{\mathrm{s}}}\] and $\mu$ is spread-out,
\begin{equation}\label{twist0}
\inf\limits_{n\in\mathbb{N}}\norm{(\nu^n\ast
\mu^n)_{\mathrm{s}}}<1 .\end{equation}  On the other hand, recall
from the proof of Lemma \ref{srfl} that the function
$\alpha\mapsto
r(\alpha):=\varrho\bigl(\hmu(\varLambda_\alpha)\bigr)$ is upper
semi-continuous and therefore attains its maximum on the compact
set $\hS$.  Fix $\xi\in\hS$ such that $r(\xi)=\max_{\alpha\in\hS}
r(\alpha)$, and observe that $\xi\neq 0$ as $0\notin\hS$.  Now
recall the argument proving \eqref{SR-5}. By Proposition
\ref{prop-1}, $\varLambda_\xi=\oplus_{i\in\mathcal{I}} U_i$ with
each $U_i$ an irreducible unitary representation of $G$. Assume
first that $\mu$ is not singular with respect to Haar measure.
Have $\epsilon>0$ with $\epsilon<1-\norm{\mu_{\mathrm{s}}}$, and
then a finite set $\mathcal{I}_1\subseteq\mathcal{I}$ such that
$\norm{\widehat{\mu_{\mathrm{a.c.}}}(U_i)}<\epsilon$ for all
$i\in\mathcal{I}_2:=\mathcal{I}\smallsetminus\mathcal{I}_1$ (use
the Riemann--Lebesgue lemma (Theorem \ref{app-thm-1})).  Then, by
\eqref{SR-5}, and as in \eqref{SR-6},
\[ \bigl[\varrho\bigl(\hmu(\varLambda_\xi)\bigr)\bigr]^n\ls
\max\limits_{i\in\mathcal{I}_1} \norm{\hmu(U_i)^n}\vee
(\epsilon+\norm{\mu_{\mathrm{s}}})^n\] for all $n\in\mathbb{N}$,
and since $\varrho\bigl(\hmu(U_i)\bigr)<1$ for each $i$, because
$\xi\neq 0$ and therefore $U_i\neq\mathbf{1}_G$ for all
$i\in\mathcal{I}$, it follows that
$\bigl[\varrho\bigl(\hmu(\varLambda_\xi)\bigr)\bigr]^n<1$ for some
sufficiently large $n$.  Hence
$\varrho\bigl(\hmu(\varLambda_\xi)\bigr)<1$, and therefore
\begin{equation}\label{twist}
\sup\limits_{\alpha\in\hS}
\Bigl[\hh(\alpha)\varrho\bigl(\hmu(\varLambda_\alpha)\bigr)\Bigr]\ls
r(\xi)=\varrho\bigl(\hmu(\varLambda_\xi)\bigr)<1 .
\end{equation} Finally, if $\mu\perp\lambda_G$, then
$\mu^n\not\perp\lambda_G$ for some $n\in\mathbb{N}$, because $\mu$
is spread-out, and the above argument with $\mu^n$ in place of
$\mu$ gives again \eqref{twist}, since also
$\varrho(x^n)=[\varrho(x)]^n$ for any element $x$ in a Banach
algebra.

It follows from \eqref{final1}, \eqref{final2}, \eqref{twist0} and
\eqref{twist}, that $f\in I_\mu$.
\end{proof}

Let $\pi_K\colon G\to K$ be the natural projection
$\pi_K(a,\kappa)=\kappa$.  Being continuous (hence Borel
measurable), $\pi_K$ induces a mapping $\pi_K\colon M(G)\to M(K)$,
given explicitly by $[\pi_K(\mu)](B)=\mu(\pi_K^{-1}(B))\
(B\in\mathcal{B}(K))$. Given a function $f\colon G\to\mathbb{C}$
and $\kappa\in K$, write $f_\kappa$ for the function
$f_\kappa\colon A\to\mathbb{C}$ given by
$f_\kappa(a)=f(a,\kappa)$; then the restriction of $\pi_K\colon
M(G)\to M(K)$ to $L^1(G)$ is given by \[\pi_K(f)(\kappa)=\int_A
f_\kappa\, d\lambda_A=\widehat{f_\kappa}(0)\qquad\text{for
$\lambda_K$-a.e.}\ \kappa\in K.\] Observe that, since
$\varLambda_0(a,\kappa)=L_K(\kappa)$ for all $a\in A$ and
$\kappa\in K$, where $L_K$ is the left regular representation of
$K$ on $L^2(K)$ (see \eqref{varLambda}), for any measure $\mu\in
M(G)$ we have that
\begin{equation}\label{piK2}\hmu(\varLambda_0)=\widehat{\pi_K(\mu)}(L_K),\end{equation} where
$\widehat{\pi_K(\mu)}(L_K)=\int_K L_K(\kappa)^\ast\,
d\pi_K(\mu)(\kappa)$ is the Fourier transform of the measure
$\pi_K(\mu)\in M(K)$ at the representation $L_K$ of $K$.  Observe
further that one also has that
\begin{equation}\label{piK3}
\pi_K(\mu\ast\nu)=\pi_K(\mu)\ast\pi_K(\nu)\end{equation} for
$\mu,\nu\in M(G)$. Indeed, \[\widehat{\pi_K(\mu\ast\nu)}(L_K)=
\widehat{(\mu\ast\nu)}(\varLambda_0)=\widehat{\nu}(\varLambda_0)\,
\widehat{\mu}(\varLambda_0) = \widehat{\pi_K(\nu)}(L_K)\,
\widehat{\pi_K(\mu)}(L_K),\] and this implies \eqref{piK3}, by
Fact \ref{GF}.   Finally observe that the kernel of the mapping
$\pi_K\colon L^1(G)\to L^1(K)$ coincides with the $L^1$-kernel of
$\varLambda_0$:
\begin{equation}\label{piK4}
\ker(\pi_K)=\ker(\varLambda_0).\end{equation} This also follows
directly from \eqref{piK2} and Fact \ref{GF}.

\begin{notation} If $U$ is any unitary representation of $G$, the
$L^1$-kernel of $U$ is $\ker([U])=\set{f\in L^1(G)\colon U(f)=0}$.
Note that, in particular for $\varLambda_0$, we also have that
$\ker(\varLambda_0)=\bigset{f\in L^1(G)\colon
\hf(\varLambda_0)=0}$.
\end{notation}

\begin{lemma}\label{lem6-1}
Let $G$ be as in Theorem \textup{\ref{main}} and let $I$ be the
ideal defined by \eqref{I}. Then $I$ is dense in the ideal
$\ker(\varLambda_0)$.
\end{lemma}

\begin{proof}
It is shown in \cite[Theorem 2]{Lud} that the hull
\[ h(\ker(\pi_K))=\bigset{\ker([U])\colon [U]\in\hG,\
\ker(\pi_K)\subseteq\ker([U])}\] of $\ker(\pi_K)$ is a set of
synthesis (see also the Remark following Lemma 3 of \cite{Lud}).
Therefore, by \eqref{piK4}, it suffices to show that the hull of
$\bar{I}$ is contained in the hull of $\ker(\pi_K)$.  Stated in
more direct terms, it suffices to show that, if $[U]\in\hG$ and
$I\subseteq\ker([U])$, then also $\ker(\pi_K)\subseteq\ker([U])$
(see \cite{Lud}, last line of the proof of Theorem 2 and Theorem
1).

Let $[U]\in\hG$ and suppose that $\ker([U])$ does not contain
$\ker(\pi_K)$.  For an irreducible unitary representation $V$ of
$K$ let $U_V$ denote the (irreducible unitary) representation of
$G$ defined by $U_V(a,\kappa):=V(\kappa)$ for all $(a,\kappa)\in
G$, and observe that
\begin{equation}\label{6.4aux}\ker(\pi_K)=\bigcap\limits_{[V]\in\hK}
\ker([U_V]),\end{equation} by \eqref{piK2} and the uniqueness of
the Fourier transform on $K$. Now $U$ is equivalent to a
sub-re\-pre\-sen\-tation of some $\varLambda_\alpha$, by
Proposition \ref{prop-1}, and by \eqref{6.4aux} and our assumption
that $\ker(\pi_K)\nsubseteq\ker([U])$ we have that $\alpha\neq 0$,
since $\varLambda_0$ is the direct sum of copies of the $U_V,\
[V]\in\hK$. There exists $f\in L^1(G)$ with $U(f)\neq 0$. Have a
neighborhood $\hW$ of $\alpha$ in $\hA$ such that the closure
$\hC$ of $\hW$ is compact and does not contain $0\in\hA$, and let
$\nu\in M(G)$ be a measure as in Lemma \ref{central} for this
$\hC$. Then $f\ast\bar{\nu}\in I\smallsetminus \ker([U])$, since
$\varLambda_\alpha(\bar{\nu})=\hnu(\varLambda_\alpha)^\ast$ and
similarly for $U$.
\end{proof}

\begin{corollary}\label{cor6-3}
Let $G$ be as in Theorem \textup{\ref{main}} and
$\ker(\varLambda_0)$ as in Lemma \textup{\ref{lem6-1}}.  If
$\mu\in M(G)$ is a spread-out probability measure with
$\varrho\bigl(\hmu(U)\bigr)< 1\ \forall\,
[U]\in\hG\smallsetminus\set{\mathbf{1}_G},$ then
$\ker(\varLambda_0)\subseteq I_\mu$.\qed
\end{corollary}

We are now ready to prove Theorem \ref{main}.

\begin{proof}[Proof of Theorem \textup{\ref{main}}]
By the Peter--Weyl theorem, $L^2(K)$ decomposes into the direct
sum
\[ L^2(K) = \bigoplus\limits_{[V]\in\hK} \mathcal{E}_{[V]},\]
where $\mathcal{E}_{[V]}$ is the finite-dimensional subspace of
$L^2(K)$ spanned by the representative functions $V_{i
j}(\kappa):=\langle V(\kappa) e_j,e_i\rangle$ of any
representation $V$ in the equivalence class $[V]$ with respect to
any basis $\set{e_1,\ldots,e_{d_{[V]}}}$ of its representation
space $\mathcal{H}_V$, and each subspace $\mathcal{E}_{[V]}$ is
invariant under the right regular representation $R_K$ of $K$ on
$L^2(K)$.  Furthermore, if $R_K^{[V]}$ denotes the
sub-representation
\[ R_K^{[V]}(\kappa):=R_K(\kappa)\vert_{\mathcal{E}_{[V]}}\] of
$R_K$, then $R_K^{[V]}$ may be decomposed into a direct sum of
$d_{[V]}$ copies from $[V]$. Next fix a Jordan normal form
decomposition of the finite-dimensional operator
$\widehat{\pi_K(\mu)}(R_K^{[V]})\colon\mathcal{E}_{[V]}\to\mathcal{E}_{[V]}$:
\[\mathcal{E}_{[V]} =
\operatorname{span}\Bigl(\phi_{1 1}^{[V]},\ldots,\phi_{1
r_1}^{[V]}\Bigr)+ \cdots+\operatorname{span}\Bigl(\phi_{p
1}^{[V]},\ldots,\phi_{p r_p}^{[V]}\Bigr)\] ($p$ and
$r_1,\ldots,r_p$ depending on $[V]$), where $\phi^{[V]}_{i j},\
i=1,\ldots,p,\, j=1,\ldots,r_i,$ are a basis of
$\mathcal{E}_{[V]}$ and satisfy
\begin{equation}\label{ge} \widehat{\pi_K(\mu)}(R_K^{[V]})\, \phi^{[V]}_{i j} =
\lambda_i^{[V]} \phi^{[V]}_{i j} + \phi^{[V]}_{i\,j-1}\quad (1\ls
i\ls r_i,\ \phi^{[V]}_{i 0}=0),\end{equation} where the
$\lambda_i^{[V]}$ are eigenvalues of the operator
$\widehat{\pi_K(\mu)}(R_K^{[V]})$. Observe that, since $R_K^{[V]}$
is a direct sum of $d_{[V]}$ copies of $[V]$, the eigenvalues
$\lambda_i^{[V]}$ are the eigenvalues of
$\widehat{\pi_K(\mu)}(V)$. Finally, observe also that, since
$\phi^{[V]}_{i j}\in\mathcal{E}_{[V]}$, we have that
\begin{equation}\label{int}
\int_K \phi^{[V]}_{i j}\, d\lambda_K = 0\end{equation} for all
$i,j$ whenever $V\neq\mathbf{1}_K$, by the Shur orthogonality
relations.

To prove the Theorem we need to show that $L^1_0(G)\subseteq
I_\mu$. By the preceding paragraph, and the Peter--Weyl theorem,
the finite linear combinations of the functions $\phi^{[V]}_{i j}$
are dense in $C(K)$ and every $L^p(K)$; hence, by standard
arguments, the functions of the form
\begin{equation}\label{dense} (a,\kappa)\mapsto g_0(a) +
\sum\limits_{[V]\in\hF} \sum\limits_{i,j} g^{[V]}_{i j}(a)\,
\phi^{[V]}_{i j}(\kappa) , \end{equation} where $\hF$ is a finite
subset of $\hK$ not containing the trivial representation
$\mathbf{1}_K$ and $g_0,g^{[V]}_{i j}\in L^1(A)$ for all
$[V]\in\hF$ and all $i,j$, are dense in $L^1(G)$. From this, it
follows that the functions of the form \eqref{dense} and with
$\int_A g_0\, d\lambda_A=0$ are dense in $L^1_0(G)$.  Therefore,
and since $I_\mu$ is closed, it suffices to show that any function
of the form \eqref{dense} with $\int_A g_0\, d\lambda_A=0$ belongs
to $I_\mu$. Since the function $f(a,\kappa):=g_0(a)$ is in
$\ker(\varLambda_0)$ if $\int_A g_0\, d\lambda_A=0$, and hence in
$I_\mu$ by Corollary \ref{cor6-3}, it is then enough to show that
any function of the form $f(a,\kappa)=g(a)\, \phi^{[V]}_{i
j}(\kappa)$ with $g\in L^1(A)$ and $V\neq\mathbf{1}_K$ belongs to
$I_\mu$. Finally, since $f\in\ker(\varLambda_0)$ if $\int_A g\,
d\lambda_A=0$, and therefore $f\in I_\mu$ by Corollary
\ref{cor6-3} again, it suffices to only consider the case $\int_A
g\, d\lambda_A\neq 0$, and then we may as well assume that $\int_A
g\, d\lambda_A=1$.

Fix $[V]\in\hK$ with $V\neq\mathbf{1}_K$, and let us suppress the
dependence of $\phi^{[V]}_{i j}$ and $\lambda_i^{[V]}$ on $[V]$
and write $\phi_{i j}$ and $\lambda_i$ instead. We will show that,
if $f_{i j}\in L^1_0(G)$ are any functions with $\pi_K(f_{i
j})=\phi_{i j}$, then $f_{i j}\in I_\mu$, which, by the above
discussion, proves the Theorem. This we will do for fixed $i$ by
induction on $j$. Fix $i$ and recall from \eqref{ge} that the
$\phi_{i j}$ satisfy $\phi_{i j}\in \mathcal{E}_{[V]}$ and
\begin{equation}\label{efns}
\widehat{\pi_K(\mu)}(R_K)\, \phi_{i j} = \lambda_i \phi_{i j} +
\phi_{i\,j-1}\quad (1\ls j\ls r_m,\ \phi_{i 0}=0),\end{equation}
where $R_K$ is the right regular representation of $K$ on
$L^2(K)$.  Obviously $f_{i 0}\in\ker(\varLambda_0)$, since
$\phi_{i 0}=0$. Assume next that $f_{i\,j-1}\in\ker(\varLambda_0)$
($j\gr 1$). Since for any measure $\nu\in M(K)$,
$\hnu(R_K)\,\phi=\phi\ast\nu$ for any $\phi\in L^2(K)$, we
conclude from \eqref{efns} and \eqref{piK3} that \[ \pi_K(f_{i
j}\ast\mu-\lambda_i f_{i j}-f_{i\,j-1})=0,\] i.e., $f_{i
j}\ast\mu-\lambda_i f_{i j}-f_{i\,j-1}\in\ker(\varLambda_0)$.
Therefore $f_{i j}\ast\mu-\lambda_i f_{i j}\in\ker(\varLambda_0)$,
because we are assuming that $f_{i\,j-1}\in\ker(\varLambda_0)$. By
Corollary \ref{cor6-3} then, $f_{i j}\ast\mu-\lambda_i f_{i j}\in
I_\mu$. Thus
\begin{equation}\label{finale}
\norm{(\lambda_i f_{i j} - f_{i j}\ast\mu)\ast\mu^n}_1 \to 0
\qquad (n\to\infty).\end{equation} Now observe that the numerical
sequence $\norm{f_{i j}\ast\mu^n}_1$ is non-increasing and bounded
by $\norm{f_{i j}}_1$.  It has therefore a limit, $a$ say, and by
\eqref{finale}, $a$ must satisfy $\abs{\lambda_i} a=a$.  Recall,
however, that $\lambda_i$ belongs to the spectrum of the operator
$\widehat{\pi_K(\mu)}(R_K^{[V]})$, which is the same as the
spectrum of $\widehat{\pi_K(\mu)}(V)=\hmu(U_V)$, where $U_V$ is
the irreducible unitary representation of $G$ defined by
$U_V(a,\kappa)=V(\kappa)$ for all $(a,\kappa)\in G$, and so by our
condition $\varrho\bigl(\hmu(U)\bigr)< 1\ \forall\,
[U]\in\hG\smallsetminus\set{\mathbf{1}_G}$ we have that
$\abs{\lambda_i}<1$.  It follows that $a=0$, and hence $f_{i j}
\in I_\mu$.  This concludes the proof of the Theorem.
\end{proof}

\section{Ergodicity by Convolutions}
\label{ergodic}

The main result of this Section is Theorem \ref{main-erg}.  Let us
begin, however, by observing that results analogous to those for
mixing of Section \ref{nec-cond} also hold for ergodicity. First,
the analogue of Proposition \ref{prop2} is the following:

\begin{proposition}
\label{prop2-ergodic} Let $G$ be a locally compact group\textup{,}
and let $\mu$ be a probability measure in $M(G)$ which is ergodic
by convolutions. Then\textup{,} for any
$[U]\in\hG\smallsetminus\set{\mathbf{1}_G},$
$n^{-1}\sum_{k=0}^{n-1}\hmu(U)^k\to 0$ in the strong operator
topology. In particular\textup{,} the number $1$ cannot be an
eigenvalue of $\hmu(U)$\textup{,} for any
$[U]\in\hG\smallsetminus\set{\mathbf{1}_G}$.
\end{proposition}

\begin{proof}
Fix $[U]\in\hG\smallsetminus\set{\mathbf{1}_G}$,
$h\in\mathcal{H}_U$, and $\epsilon>0$. As in the proof of \1\ of
Proposition \ref{prop2}, there exist $g\in L^1_0(G)$ and
$h'\in\mathcal{H}_U$ for which $\norm{\hg(U) h'-h} < \epsilon/2$.
It follows that
\begin{align*}
\Biggnorm{\frac{1}{n}\sum\limits_{k=0}^{n-1} \hmu(U)^k h} &\ls
\Biggnorm{\frac{1}{n}\sum\limits_{k=0}^{n-1} \hmu(U)^k \hg(U) h'}
+ \frac{1}{n}\sum\limits_{k=0}^{n-1}\bignorm{\hmu(U)^k (\hg(U) h'
- h)}
\\ &< \Biggnorm{\frac{1}{n}\sum\limits_{k=0}^{n-1} g\ast\mu^n}_1 +
\epsilon/2 ,
\end{align*}
and this is $<\epsilon$ for all sufficiently large $n$, since
$g\in L^1_0(G)$ and $\mu$ is ergodic by convolutions.
\end{proof}

Again, one can say more about spread-out measures. Combining
Proposition \ref{prop2-ergodic} with Lemma \ref{Yosida-Kakutani-1}
and Lemma \ref{lem1}, one obtains the following:

\begin{corollary}\label{cor-erg}
Let $G$ be a locally compact \textup{CCR} group\textup{,} and let
$\mu$ be a spread-out probability measure in $M(G)$ which is
ergodic by convolutions. Then $1\notin\sigma\bigl(\hmu(U)\bigr)$
for any $[U]\in\hG\smallsetminus\set{\mathbf{1}_G}$.
\end{corollary}

Our main result concerning ergodicity by convolutions is the
analogue of Theorem \ref{main} for ergodicity, namely that the
necessary condition of the above Corollary is also sufficient for
ergocity in motion groups.

\begin{theorem}\label{main-erg}
Let $G=A\times_{\varphi} K$ be a motion group with $G$ acting
regularly on $\hA$.  Then\textup{,} if $\mu\in M(G)$ is a
spread-out probability measure\textup{,} $\mu$ is ergodic by
convolutions if and only if $1\notin\sigma\bigl(\hmu(U)\bigr)$ for
any $[U]\in\hG\smallsetminus\set{\mathbf{1}_G}$.
\end{theorem}

To prove the Theorem we only have to show that, for a spread-out
probability measure $\mu\in M(G)$,
$1\notin\bigcup_{[U]\in\hG\smallsetminus\set{\mathbf{1}_G}}\sigma\bigl(\hmu(U)\bigr)$
implies ergodicity, the other direction being a consequence of
Corollary \ref{cor-erg}. Let \[ J_\mu:=\Biggset{f\in L^1(G)\colon
\lim\limits_{n\to\infty}\Biggnorm{\frac{1}{n}\sum\limits_{k=0}^{n-1}
f\ast\mu^k}_1=0}.\]  Then $J_\mu$ is a closed left ideal in
$L^1(G)$ contained in $L^1_0(G)$, and we have to show that
$J_\mu=L^1_0(G)$.  As in the case of mixing, we will first show
that $\ker(\varLambda_0)\subseteq J_\mu$ and then deduce from this
that all of $L^1_0(G)$ is contained in $J_\mu$.  To prove the
first assertion, we will use the following Lemma, whose proof we
postpone to the end of this Section.

\begin{lemma}\label{ErgLem1}
Let $\mu$ be a probability measure in $M(G)$, with
$\mu\not\perp\lambda_G$ and for which
$1\notin\sigma\bigl(\hmu(U)\bigr)\ \forall\,
[U]\in\hG\smallsetminus\mathbf{1}_G$.  Let also $\nu\in Z(M(G))$
be a measure as in Lemma \textup{\ref{central},} for some
$\hC\in\widehat{\mathcal{C}}_0$.  Then $\delta_e-\mu\ast\nu$ is
invertible in $D^\ast(G)$, where $\delta_e$ is the Dirac
point-mass at the neutral element $e$ of $G$.
\end{lemma}

\begin{lemma}\label{ErgLem2}
Let $G$ be as in Theorem \textup{\ref{main-erg}} and
$\ker(\varLambda_0)$ as in Lemma \textup{\ref{lem6-1}}.  If
$\mu\in M(G)$ is a spread-out probability measure for which
$1\notin\sigma\bigl(\hmu(U)\bigr)\ \forall\,
[U]\in\hG\smallsetminus\mathbf{1}_G,$ then
$\ker(\varLambda_0)\subseteq J_\mu$.
\end{lemma}

\begin{proof}
Since $J_\mu$ is closed, it suffices, by Lemma \ref{lem6-1}, to
show that $I\subseteq J_\mu$, where $I$ is defined in \eqref{I}.
Assume first that $\mu$ is not singular with respect to Haar
measure $\lambda_G$. Next observe that
\[
\Biggnorm{\frac{1}{n}\sum\limits_{k=0}^{n-1}(f-f\ast\mu)\ast\mu^k}_1
= \frac{1}{n}\norm{f-f\ast\mu^n}_1\ls \frac{2}{n}\,\norm{f}_1
\qquad\forall\, f\in L^1(G),\] and therefore it suffices to show
that, for each $g\in I$, there exists $f\in L^1(G)$ such that
$g=f-f\ast\mu$.

Fix $g\in I$.  Then there exists a compact set $\hC\in
\widehat{\mathcal{C}}_0$ such that $\hg(\varLambda_\alpha)=0$ for
$\alpha$ not in $\hC$.  Fix a measure $\nu\in Z(M(G))$ as in Lemma
\ref{central} for this $\hC$, and set $\mu^\prime:=\mu\ast\nu$.
Since $\mu\not\perp\lambda_G$ we also have that
$\mu^\prime\not\perp\lambda_G$.  Now recall the argument in the
proof of Lemma \ref{srfl}.  Write
$\mu^\prime=\mu^\prime_{\mathrm{a.c.}}+\mu^\prime_{\mathrm{s}}$
and observe that $\delta_e-\mu^\prime_{\mathrm{s}}$ is an
invertible element of $M(G)$, because
$\norm{\mu^\prime_{\mathrm{s}}}<1$.  By the preceding lemma we
also have that $\delta_e-\mu^\prime$ is invertible in $D^\ast(G)$.
Set $\nu^\prime :=
(\delta_e-\mu^\prime_{\mathrm{s}})^{-1}\ast\mu^\prime_{\mathrm{a.c.}}$;
since
\begin{equation}\label{A}
(\delta_e-\mu^\prime_{\mathrm{s}})^{-1}\ast (\delta_e-\mu^\prime)
= \delta_e-\nu^\prime \end{equation} and the left side is
invertible in $D^\ast(G)$, $\delta_e-\nu^\prime$ is invertible in
$D^\ast(G)$.  But $\nu^\prime\in L^1(G)$, and therefore
$\delta_e-\nu^\prime$ is invertible in $D^\ast(G)$ if and only if
it is invertible in $M(G)$, by Proposition \ref{spectra}.
Therefore $\delta_e-\nu^\prime\in M(G)$. It follows from this and
\eqref{A} that
\[(\delta_e-\mu^\prime)^{-1} = (\delta_e-\nu^\prime)^{-1}\ast
(\delta_e-\mu^\prime_{\mathrm{s}})^{-1} \in M(G)\ast M(G)=M(G).\]
Now set $f:=g\ast (\delta_e-\mu^\prime)^{-1}\in L^1(G)$.  Then
$f-f\ast\mu^\prime=g$.  Furthermore, since
$\hf(\varLambda_\alpha)=\Bigl[\widehat{\delta_e}(\varLambda_\alpha)
-\widehat{(\mu^\prime)}(\varLambda_\alpha)\Bigr]^{-1}
\hg(\varLambda_\alpha)=0$ for $\alpha\in\hA\smallsetminus\hC$ and
$\hnu(\varLambda_\alpha)=I_{L^2(K)}$ on $\hC$, it follows that
$f\ast\nu=f$, and therefore $f-f\ast\mu^\prime=f-f\ast\mu$.

If $\mu$ is singular with respect to Haar measure, then some power
$\mu^m$ of $\mu$ is not singular, because $\mu$ is assumed to be
spread-out. Replacing $\mu$ by $\mu^m$ in the above argument
yields, for a given $g\in I$, a function $f\in L^1(G)$ for which
$g=f-f\ast\mu^m$.  But then
\begin{align*}
\Biggnorm{\frac{1}{n}\sum\limits_{k=0}^{n-1}g\ast\mu^k}_1 &=
\Biggnorm{\frac{1}{n}\sum\limits_{k=0}^{n-1}(f-f\ast\mu^m)\ast\mu^k}_1
\\ &= \frac{1}{n}\Biggnorm{\sum\limits_{k=0}^{m-1}f\ast\mu^k-
\sum\limits_{k=n}^{n+m-1}f\ast\mu^k}_1 \ls
\frac{2m}{n}\,\norm{f}_1 ,\end{align*} and therefore $g\in J_\mu$
again.
\end{proof}

\begin{proof}[Proof of Theorem \textup{\ref{main-erg}}]
By the argument in the proof of Theorem \ref{main} it suffices to
show that, if $f_{i j}\in L^1_0(G)$ are any functions with
$\pi_K(f_{i j})=\phi_{i j}^{[V]}$ for some
$V\in\hK\smallsetminus\set{\mathbf{1}_K}$, then $f_{i j}\in
J_\mu$, where the notation is as in that proof.  Let us suppress
the dependence on $[V]$ and write $\lambda_i$ and $\phi_{i j}$
instead of $\lambda_i^{[V]}$ and $\phi_{i j}^{[V]}$ again.  As
shown in the proof of Theorem \ref{main}, $f_{i
j}\ast\mu-\lambda_i f_{i j}$ is in $\ker(\varLambda_0)$, for all
$j\in\set{1,\ldots,r_i}$, so by Lemma \ref{ErgLem2}
\[ \frac{1}{n}\Biggnorm{\sum\limits_{k=0}^{n-1}
(\lambda_i f_{i j}-f_{i j}\ast\mu)\ast\mu^k}_1\to 0
\qquad(n\to\infty),\] and hence
\[ \frac{1}{n}\Biggnorm{(\lambda_i-1)\sum\limits_{k=1}^{n-1}
f_{i j}\ast\mu^k+\lambda_i f_{i j}-f_{i j}\ast\mu^n}_1\to 0
\qquad(n\to\infty),\] for all $j$.  From this, it follows that
\[\frac{1}{n}\,\abs{\lambda_i-1}\,\Biggnorm{\sum\limits_{k=1}^{n-1}
f_{i j}\ast\mu^k}_1\to 0 \qquad(n\to\infty),\] and as
$\lambda_i\neq 1$, because we are assuming that
$1\notin\sigma\bigl(\hmu(U)\bigr)\ \forall\,
[U]\in\hG\smallsetminus\set{\mathbf{1}_G}$ and that
$V\neq\mathbf{1}_K$, we must have that $f_{i j}\in J_\mu$.
\end{proof}

\begin{proof}[Proof of Lemma \textup{\ref{ErgLem1}}]
Fix a complete set of mutually inequivalent, irreducible, unitary
representations of $G$, denote it by $\mathcal{U}_G$, and consider
the unital $C^\ast$-algebra
\[\boldsymbol{C}(\hG):=\biggset{(T_U)_{U\in\mathcal{U}_G}\colon
T_U\in \boldsymbol{B}(\mathcal{H}_U),\
\sup_{U\in\mathcal{U}_G}\norm{T_U} < \infty} ,\] with norm
$\sup_{U\in\mathcal{U}_G}\norm{T_U}$, and pointwise operations.
Since $\iota\colon D^\ast(G)\to \boldsymbol{C}(\hG)$ given by
$\iota(x):= (U(x))_{U\in\mathcal{U}_G}$ is an isometry,
$\delta_e-\mu\ast\nu$ is invertible in $D^\ast(G)$ if and only if
$(U(\delta_e-\mu\ast\nu))_{U\in\mathcal{U}_G}$ is invertible in
$\boldsymbol{C}(\hG)$ \cite[Proposition 1.23]{Fol}.

Recall Proposition \ref{prop-1} and write
$[U]\in\langle\varLambda_\alpha\rangle$ if $U$ is unitarily
equivalent to a sub-representation of $\varLambda_\alpha$. Let
$\hS$ denote the support of $\hh$, which, recall, is compact and
does not contain $0\in\hA$. Then
$\varLambda_\alpha(\delta_e-\mu\ast\nu)=I_{L^2(K)} -
\hh(\alpha)\varLambda_\alpha(\mu)=I_{L^2(K)}$ for $\alpha$ outside
$\hS$, and therefore
\begin{equation}\label{aa}[U(\delta_e-\mu\ast\nu)]^{-1}=I_{\mathcal{H}_U}\qquad \text{for}\
[U]\in\langle\varLambda_\alpha\rangle\ \text{with}\
\alpha\in\hA\smallsetminus\hS.\end{equation}

Next let $\boldsymbol{B}^\prime(L^2(K))$ denote the invertible
elements of $\boldsymbol{B}(L^2(K))$, and recall that the mapping
$x\mapsto x^{-1}$ is continuous on $\boldsymbol{B}^\prime(L^2(K))$
\cite[Theorem 1.4]{Fol}.  Since the mapping
$\jmath\colon\hA\to\boldsymbol{B}(L^2(K))$ with
$\jmath(\alpha):=\varLambda_\alpha(\delta_e-\mu\ast\nu)$ is also
continuous (Theorem \ref{app-thm-2}), if we show that
$\jmath(\hA)\subseteq\boldsymbol{B}^\prime(L^2(K))$, then
$\alpha\mapsto
\bignorm{[\varLambda_\alpha(\delta_e-\mu\ast\nu)]^{-1}}$ will be
continuous on $\hA$, and hence bounded on the compact $\hS$; thus
we will have that
\[\bignorm{[U(\delta_e-\mu\ast\nu)]^{-1}} \ls
\sup\limits_{\alpha\in\hA}
\bignorm{[\varLambda_\alpha(\delta_e-\mu\ast\nu)]^{-1}}<\infty
\qquad \text{for}\ [U]\in\langle\varLambda_\alpha\rangle\
\text{and}\ \alpha\in\hS,\] and this together with \eqref{aa} will
show that
$([U(\delta_e-\mu\ast\nu)]^{-1})_{U\in\mathcal{U}_G}\in\boldsymbol{C}(\hG).$

It remains to show that
$\jmath(\hA)\subseteq\boldsymbol{B}^\prime(L^2(K))$, i.e., that
$\jmath(\alpha)=\varLambda_\alpha(\delta_e-\mu\ast\nu)$ is
invertible for each $\alpha$, and by the line preceding \eqref{aa}
it suffices to only consider $\alpha\neq 0\in\hA$.  Fix such an
$\alpha$.  Then for each $[U]\in\langle\varLambda_\alpha\rangle$
we have that $1\notin\sigma(U(\mu\ast\nu))$; for if
$0\ls\hh(\alpha)<1$ then
$\norm{U(\mu\ast\nu)}=\hh(\alpha)\norm{U(\mu)}<1$, and if
$\hh(\alpha)=1$ then $U(\mu\ast\nu)=U(\mu)$ and
$1\notin\sigma(U(\mu))$ by hypothesis, since $\alpha\neq 0$
implies that $U\neq\mathbf{1}_G$. Thus $U(\delta_e-\mu\ast\nu)$ is
invertible for each $[U]\in\langle\varLambda_\alpha\rangle$. By
Proposition \ref{prop-1},
$\varLambda_\alpha=\oplus_{i\in\mathcal{I}} U_i$ with each $U_i$
an irreducible unitary representation of $G$. Recall also that we
are assuming that $\mu$ is not singular with respect to Haar
measure. Have $\epsilon>0$ with
$\epsilon<1-\norm{\mu_{\mathrm{s}}}$, and then a finite set
$\mathcal{I}_1\subseteq\mathcal{I}$ such that
$\norm{U_i(\mu_{\mathrm{a.c.}})}<\epsilon$ for all
$i\in\mathcal{I}_2:=\mathcal{I}\smallsetminus\mathcal{I}_1$ (use
the Riemann--Lebesgue lemma (Theorem \ref{app-thm-1})).  Then
\[\norm{U_i(\mu\ast\nu)}=\hh(\alpha) \norm{U_i(\mu)}\ls
\norm{U_i(\mu)}\ls
\norm{U_i(\mu_{\mathrm{a.c.}})}+\norm{U_i(\mu_{\mathrm{s}})} \ls
\epsilon + \norm{\mu_{\mathrm{s}}}
\] for each $i\in\mathcal{I}_2$, and therefore, if we set
$U:=\bigoplus_{i\in\mathcal{I}_2} U_i$, then
\[ \norm{U(\mu\ast\nu)}\ls \epsilon + \norm{\mu_{\mathrm{s}}}<1.\]
It follows that $U(\delta_e-\mu\ast\nu)$ is invertible on
$\mathcal{H}_U$, and since
\[\varLambda_\alpha(\delta_e-\mu\ast\nu)=\left[\bigoplus_{i\in\mathcal{I}_1}
U_i(\delta_e-\mu\ast\nu)\right] \oplus\, U(\delta_e-\mu\ast\nu)\]
is a finite sum with each summand invertible, it follows that
$\varLambda_\alpha(\delta_e-\mu\ast\nu)$ is invertible.
\end{proof}

\section{Weak Mixing}
\label{weakmix}

In \cite{Rblatt} Rosenblatt observes that, by the work of Foguel
\cite{Fog}, weak mixing is actually equivalent to mixing in
Abelian groups, and asks whether this remains true for more
general groups. The answer turns out in the affirmative for
spread-out measures on motion groups with $G$ acting regularly on
$\hA$. To prove this, it suffices to show the following result and
then refer to Theorem \ref{main}:

\begin{proposition}\label{wmix}
Let $G$ be a locally compact \textup{CCR}  group\textup{,} and let
$\mu$ be a spread-out probability measure in $M(G)$ which is
weakly mixing by convolutions. Then $\varrho\bigl(\hmu(U)\bigr)<1$
for any $[U]\in\hG\smallsetminus\set{\mathbf{1}_G}$.
\end{proposition}

\begin{proof}
By Aaronson et al. \cite{Aetal}, weak mixing is equivalent to the
following condition: $\mu\ast h=\lambda h$ with $h\in L^\infty(G)$
and $\abs{\lambda}=1$ implies that $\lambda=1$ and $h$ is constant
$\lambda_G$-a.e.  Assume that
$\lambda\in\sigma\bigl(\hmu(U)\bigr)$ for some
$[U]\in\hG\smallsetminus\set{\mathbf{1}_G}$ and
$\lambda\in\mathbb{C}$ with $\abs{\lambda}=1$.  Then also
$\overline{\lambda}\in \sigma(U(\mu))$. Since $\mu$ is spread-out,
and hence $U(\mu)$ is quasi-compact, $\overline{\lambda}$ must be
an eigenvalue of $U(\mu)$ (see Lemma \ref{Yosida-Kakutani-1}). Let
$u\in\mathcal{H}_U$ be an eigenvector for $\overline{\lambda}$
with $\norm{u}=1$, i.e., assume that
$U(\mu)u=\overline{\lambda}u$, and set $h(x):=\langle
U(x)u,u\rangle\ (x\in G)$.  Then $\mu\ast h = \lambda h$, so if
$\mu$ is weakly mixing we must have that $\lambda=1$ and that $h$
is constant (since it is also continuous).  But if $U\neq
\mathbf{1}_G$, $h(x)=\langle U(x)u,u\rangle$ can not be constant.
\end{proof}

\begin{corollary}\label{mix=wmix}
Let $G=A\times_{\varphi} K$ be a motion group with $G$ acting
regularly on $\hA$.  Then\textup{,} if $\mu\in M(G)$ is a
spread-out probability measure\textup{,} $\mu$ is mixing by
convolutions if and only if it is weakly mixing by convolutions.
\end{corollary}

\section{Final Remarks}
\label{cr}

A probability measure $\mu\in M(G)$ on a locally compact group $G$
is \emph{adapted\/} if it is not concentrated on a closed proper
subgroup of $G$, and \emph{strictly aperiodic\/} if it is not
concentrated on a coset of a normal, closed, proper subgroup of
$G$.

Consider the following conditions for a probability measure
$\mu\in M(G)$ on a locally compact group $G$:
\begin{itemize}
\item[(E)]
$\mu$ is ergodic.
\item[(M)]
$\mu$ is mixing.
\item[(A)]
$\mu$ is adapted.
\item[(ASA)]
$\mu$ is aperiodic, i.e., adapted and strictly aperiodic.
\item[(S)]
$1\notin\sigma\bigl(\hmu(U)\bigr)\ \forall\,
[U]\in\hG\smallsetminus\set{\mathbf{1}_G}$.
\item[(SR)]
$\varrho\bigl(\hmu(U)\bigr)< 1\ \forall\,
[U]\in\hG\smallsetminus\set{\mathbf{1}_G}$.
\end{itemize}

In Abelian and compact groups, it is known that
$\text{(E)}\Leftrightarrow\text{(A)}\Leftrightarrow\text{(S)}$,
and that
$\text{(M)}\Leftrightarrow\text{(ASA)}\Leftrightarrow\text{(SR)}$
(see \cite{ChD}, \cite{Fog}, \cite[Theorem 2 and Remark 1]{Ram}
and \cite[Proposition 1.2 and Theorem 1.4]{Rblatt} for Abelian
groups, and \cite{Kaw-Ito}, \cite{Strom}, \cite[Theorem
V.5.2]{MRos} and  \cite[2.5.14]{Heyer} for compact groups). The
equivalence $\text{(E)}\Leftrightarrow\text{(A)}$ is also known
for spread-out measures in locally compact, compactly generated,
second countable groups of polynomial growth (\cite{Jaw1}).  In
fact, adaptedness is necessary for ergodicity, and aperiodicity
necessary for mixing, in \emph{any} locally compact group (see
\cite[p.\ 33 and p.\ 38]{Rblatt}).   On the other hand, the
example on p.~40 of \cite{Rblatt} shows that these conditions are
no longer sufficient for ergocity and mixing, respectively, in
arbitrary groups. Here, we make the observation that, for this
example, there are in fact non-trivial irreducible unitary
representations of the underlying group for which condition (S)
fails.

\begin{example*}[Rosenblatt \cite{Rblatt}]
The underlying group in this example is the semi-direct product
$G=\mathbb{Z}^2\times_\varphi\mathbb{Z}$, where $\mathbb{Z}$ acts
on $\mathbb{Z}^2$ through the automorphisms
$\varphi_k(n_1,n_2)=(n_1,n_2)\varGamma$, where
$\varGamma=\begin{pmatrix} 1 & 2 \\ 2 & 3\end{pmatrix}$.  Set
$a=(0,0,1)$, $b=(1,2,1)$, and $ c=(2,3,1)$, and consider the
probability measure
$\mu:=\tfrac{1}{4}(\delta_a+\delta_b+\delta_{bc}+\delta_{c^2})$.
It is shown in \cite{Rblatt} that $\mu$ is aperiodic, yet neither
mixing nor ergodic.  Note that $\mu$ is certainly spread-out, as
$G$ is discrete.  For
$(t_1,t_2)\in\mathbb{T}^2=\mathbb{R}^2/\mathbb{Z}^2=\widehat{\mathbb{Z}^2}$,
let $\varLambda_{(t_1,t_2)}$ be the representation of $G$ on
$\ell^2(\mathbb{Z})$ given by
\[ [\varLambda_{(t_1,t_2)}(n_1,n_2,k)\phi](m) =
\exp(2\pi i (t_1,t_2)\varGamma^{-m} (n_1,n_2)^\prime)\,\phi(m-k)\]
with the dash denoting transpose, i.e.,
$(n_1,n_2)^\prime=\begin{pmatrix}n_1 \\ n_2\end{pmatrix}$. The
representation $\varLambda_{(t_1,t_2)}$ is easily seen to be
unitary, and it is also irreducible, as can be seen using Shur's
lemma. Furthermore, direct computation shows that
\begin{multline*} [\varLambda_{(t_1,t_2)}(\mu) \phi] (m) =
\tfrac{1}{4} [1+\exp(2\pi i
(t_1,t_2)\varGamma^{-m}(1,2)^\prime)]\, \phi(m-1)\\ + \tfrac{1}{4}
[\exp(2\pi i (t_1,t_2)\varGamma^{-m}(9,15)^\prime) + \exp(2\pi i
(t_1,t_2)\varGamma^{-m}(10,16)^\prime)]\, \phi(m-2)
\end{multline*} $(\phi\in\ell^2(\mathbb{Z}),\ m\in\mathbb{Z})$.
So, if $\phi_n\in\ell^2(\mathbb{Z})$ is the function with
\[ \phi_n(m)=\begin{cases}1/\sqrt{n}, & 0\ls m <n\\ 0, & \text{otherwise},\end{cases}\]
then $\norm{\phi_n}=1$, and, by direct computation again,
\begin{align*}&\norm{\varLambda_{(t_1,t_2)}(\mu) \phi_n-\phi_n}^2 =\frac{1}{n} +
\frac{1}{16n} \abs{1+\exp(2\pi i
(t_1,t_2)\varGamma^{-1}(1,2)^\prime)-4}^2\\ &\qquad +
\frac{1}{16n} \abs{\exp(2\pi i
(t_1,t_2)\varGamma^{-n-1}(9,15)^\prime) + \exp(2\pi i
(t_1,t_2)\varGamma^{-n-1}(10,16)^\prime)}^2\\ &\qquad +
\frac{1}{16n} \abs{1 + \exp(2\pi i
(t_1,t_2)\varGamma^{-n}(1,2)^\prime)\\ &\qquad\quad + \exp(2\pi i
(t_1,t_2)\varGamma^{-n}(9,15)^\prime) + \exp(2\pi i
(t_1,t_2)\varGamma^{-n}(10,16)^\prime)}^2\\ &\qquad +
\frac{1}{16n} \sum\limits_{j=2}^{n-1} \abs{1 + \exp(2\pi i
(t_1,t_2)\varGamma^{-n}(1,2)^\prime)\\ &\qquad\quad + \exp(2\pi i
(t_1,t_2)\varGamma^{-n}(9,15)^\prime) + \exp(2\pi i
(t_1,t_2)\varGamma^{-n}(10,16)^\prime)-4}^2.
\end{align*}
Now observe that $\lambda=\sqrt{5}-2$ is an eigenvalue of
$\varGamma^{-1}$, and choose $(t_1,t_2)$ to be any eigenvector
corresponding to $\lambda$.  Then, since $\abs{\lambda}<1$,
\[ \exp(2\pi i
(t_1,t_2)\varGamma^{-n}(n_1,n_2)^\prime) = \exp(2\pi i
\lambda^n(t_1,t_2) (n_1,n_2)^\prime)\to 1\] as $n\to\infty$ for
any $(n_1,n_2)\in\mathbb{Z}^2$, and hence
$\norm{\varLambda_{(t_1,t_2)}(\mu) \phi_n-\phi_n}^2\to 0$ as
$n\to\infty$.  Since $\norm{\phi_n}=1$, this shows that
$1\in\sigma\bigl(\varLambda_{(t_1,t_2)}(\mu)\bigr)$.
\end{example*}

Concluding, let us also mention the following facts in relation to
the above: \textit{in a locally compact \textup{CCR}
group\textup{,} strict aperiodicity is equivalent to the condition
\textup{(SR)} for adapted spread-out measures}. This may be proved
along the lines of Lemma 4.3 of \cite{AG}, using also Lemma
\ref{Yosida-Kakutani-1} and Lemma \ref{lem1} of the present paper,
and the fact that, in any locally compact group $G$, if $H$ is a
closed, normal, proper subgroup of $G$, then there exists a
non-trivial irreducible unitary representation $U$ of $G$ which is
identically equal to the identity operator when restricted to the
subgroup $H$. One can also prove along the lines of \cite[Lemma
4.3]{AG} that: \textit{in a \textup{CCR} group\textup{,}
$\textup{(A)}\Rightarrow\textup{(S)}$ for spread-out measures\/}.
Hence also, the equivalences
$\text{(E)}\Leftrightarrow\text{(A)}\Leftrightarrow\text{(S)}$ and
$\text{(M)}\Leftrightarrow\text{(ASA)}\Leftrightarrow\text{(SR)}$
also hold for spread-out measures in motion groups like the ones
studied in this paper.  In fact there is a direct argument showing
that $\text{(S)}\Rightarrow\text{(A)}$ for spread-out measures on
motion groups, which we now present; this, when combined with
Jaworski's theorem \cite[Corollary 3.8]{Jaw1} and the
aforementioned implication $\textup{(A)}\Rightarrow\textup{(S)}$
for spread-out measures in CCR groups, yields another proof of
Theorem \ref{main-erg}, when the motion group in question is
second countable (see also Remark 3.9 of \cite{Jaw1}).

\begin{lemma}\label{thelemma}
Let $G=A\times_{\varphi} K$ be a motion group with $G$ acting
regularly on $\hA$\textup{,} and let $H$ be an open proper
subgroup of $G$.  Then there exists
$[U]\in\hG\smallsetminus\set{\mathbf{1}_G}$ such that the
restriction of $U$ to $H$ has a fixed vector\textup{:} $\exists\,
u\in\mathcal{H}_U\ \textup{s.t.}\ U(x)u=u\ \forall\, x\in H$.
\end{lemma}

\begin{proof}
For notational simplicity we shall identify $A$ and
$A\times\set{1}$ in what follows. Assume first that $H A\neq G$.
Since $H$ is open, $HA$ is open, hence its projection
$\pi_K(HA)=\set{\kappa\in K\colon (a,\kappa)\in HA\ \text{for
some}\ a\in A}$ onto $K$ is an open subgroup of $K$; thus
$\pi_K(HA)$ is also a closed subgroup of $K$.  By the Frobenious
reciprocity theorem, there exists a non-trivial, irreducible
unitary representation $[V]\in\hK$ which has a fixed vector when
restricted to $\pi_K(HA)$.  Then $U(a,\kappa):=V(\kappa)$ for all
$a\in A$ and $\kappa\in K$ is the sought for representation of
$G$.

Next assume that $HA=G$. Then $H\cap A$ is normal in $G$.  Set
$\tilde{G}:=G/H\cap A$, $\tilde{A}:=A/H\cap A$, and
$\tilde{H}:=H/H\cap A$.  Note that $\tilde{A}$ and $\tilde{H}$ are
closed subgroups of $\tilde{G}$, and that $\tilde{A}$ is normal in
$\tilde{G}$ and $\tilde{H}\cap\tilde{A}$ is trivial.  It follows
that $\tilde{G}$ is the semi-direct product
$\tilde{G}=\tilde{A}\ltimes\tilde{H}$, in fact with the factor
$\tilde{H}$ compact, since $H$ is open and $\tilde{H}=H/H\cap
A\simeq HA/A =G/A$.  Now let $\tilde{\alpha}$ be a non-trivial
character of $\tilde{A}$, and let
$\tilde{\varLambda}_{\tilde{\alpha}}=\operatorname{ind}^{\tilde{G}}_{\tilde{A}}
(\tilde{\alpha})$ be the representation of $\tilde{G}$ on
$L^2(\tilde{H})$ obtained as in \eqref{varLambda}.  It is readily
seen that the constant function $\mathbf{1}_{\tilde{H}}$ is a
fixed vector for the restriction of
$\tilde{\varLambda}_{\tilde{\alpha}}$ on $\tilde{H}$:
$\tilde{\varLambda}_{\tilde{\alpha}}(x) \mathbf{1}_{\tilde{H}} =
\mathbf{1}_{\tilde{H}}$ for all $x\in\tilde{H}$.  It follows from
Proposition \ref{prop-1} that $\tilde{U}(x)\phi = \phi\ \forall\,
x\in\tilde{H}$ for some non-zero $\phi\in L^2(\tilde{H})$, for
some irreducible sub-representation of
$\tilde{\varLambda}_{\tilde{\alpha}}$, and since $\tilde{\alpha}$
is non-trivial $\tilde{U}$ is non-trivial, again by Proposition
\ref{prop-1} (note that $\tilde{G}$ acts regularly on the
characters of $\tilde{A}$ because $G$ acts regularly on $\hA$).
Now lift $\tilde{U}$ to $G$: the representation
$U(a,\kappa):=\tilde{U}((a,\kappa)(H\cap A))$ is the desired
representation of $G$.
\end{proof}

To obtain the asserted implication
$\text{(S)}\Rightarrow\text{(A)}$ for a spread-out measure $\mu$
on a motion group $G=A\times_{\varphi} K$ with $G$ acting
regularly on $\hA$, we then argue as follows.  It is easy to see
that, when $\mu$ is spread-out, some power $\mu^n$ of $\mu$ must
dominate a positive multiple of Haar measure $\lambda_G$ on some
open set; from this it follows that the smallest closed subgroup
$H$ of $G$ with $\mu(H)=1$ is also open.  Thus if $\mu$ is not
adapted, there exists a non-trivial irreducible unitary
representation $U$ of $G$ which has a non-trivial fixed vector $u$
say when restricted to $H$. It follows that $\hmu(U)u=u$, since
$\mu(H)=1$, and thus $1\in\sigma\bigl(\hmu(U)\bigr)$.

\appendix\section*{Appendix}
\setcounter{section}{1}

In this Appendix, we give the proofs of Theorems \ref{app-thm-1}
and \ref{app-thm-2}.

\begin{proof}[Proof of Theorem \textup{\ref{app-thm-1}}]  \1\  First, since
$\hf(\varLambda_\alpha)=\varLambda_\alpha(\bar{f})^\ast$ for any
$f\in L^1(G)$, we may as well consider $\varLambda_\alpha(f)$
instead of $\hf(\varLambda_\alpha)$.  Second, since
$\varLambda_\alpha$ is a unitary representation, and hence
\[
\norm{\varLambda_\alpha(f)-\varLambda_\alpha(g)} \ls \norm{f-g}_1
\qquad (f,g\in L^1(G),\ \alpha\in\hA) ,\] it suffices to only
consider functions of the form
\begin{equation}\label{app-1} g(a,\kappa)=\sum\limits_{V\in\mathfrak{F}}
\sum\limits_{i,j=1}^{d_{[V]}} g^V_{i j}(a)\, V_{i j}(\kappa) ,
\end{equation} where $\mathfrak{F}$ is a finite set of mutually
inequivalent, irreducible, unitary representations of $K$, the
$V_{i j}$ are the representative functions of a $V\in\mathfrak{F}$
with respect to some basis of $\mathcal{H}_V$, and $g^V_{i j}\in
L^1(A)$ for all $V\in\mathfrak{F}$ and $1\ls i,j\ls d_{[V]}$,
because the functions of the form \eqref{app-1} are dense in
$L^1(G)$. Fix a $g$ as in \eqref{app-1}.

For an arbitrary function $h\in L^1(G)$ write $h_\kappa$ for the
function $h_\kappa(a):=h(a,\kappa)$ on $A$; then $h_\kappa\in
L^1(A)$ for $\lambda_K$-a.e.\ $\kappa$ in $K$, and for $\phi\in
L^2(K)$,
\begin{align}\label{app-1.5}
[\varLambda_\alpha(h)\phi](\tkappa) &= \int_K \int_A
h(a,\kappa)\cdot \langle a,\varphi_{\tkappa}(\alpha) \rangle\cdot
\phi(\kappa^{-1}\tkappa)\, d\lambda_A(a)\, d\lambda_K(\kappa)\\ &=
\int_K \widehat{h_\kappa}(-\varphi_{\tkappa}(\alpha))\cdot
\phi(\kappa^{-1}\tkappa)\, d\lambda_K(\kappa) ,
\nonumber\end{align} where
$\widehat{h_\kappa}(\varphi_{\tkappa}(\alpha))$ is the Fourier
transform of the function $h_\kappa\colon A\to \mathbb{C}$ on $A$
at the character $\varphi_{\tkappa}(\alpha)$ of $A$. Hence
\begin{align*} \norm{\varLambda_\alpha(h) \phi}^2_{L^2(K)} &=
\int_K \biggabs{\int_K
\widehat{h_\kappa}(-\varphi_{\tkappa}(\alpha))\,
\phi(\kappa^{-1}\tkappa)\, d\lambda_K(\kappa)}^2\,
d\lambda_K(\tkappa)
\\ &\ls \norm{\phi}^2_{L^2(K)}\, \int_K \int_K
\Bigabs{\widehat{h_\kappa}(-\varphi_{\tkappa}(\alpha))}^2
d\lambda_K(\tkappa)\, d\lambda_K(\kappa) ,\end{align*} and
therefore
\begin{equation}\label{app-2}
\norm{\varLambda_\alpha(h)} \ls \left( \int_K \int_K
\Bigabs{\widehat{h_\kappa}(-\varphi_{\tkappa}(\alpha))}^2
d\lambda_K(\tkappa)\, d\lambda_K(\kappa) \right)^{1/2}
.\end{equation}  For the function $g$ defined by \eqref{app-1},
\begin{equation}\label{app-2.5} \widehat{g_{\kappa}}(\alpha) =
\sum\limits_{V\in\mathfrak{F}} \sum\limits_{i,j=1}^{d_{[V]}}
\widehat{(g^V_{i j})}(\alpha)\, V_{i j}(\kappa) ,\end{equation}
whence by \eqref{app-2}
\begin{equation}\label{app-3} \norm{\varLambda_\alpha(g)}^2 \ls
\sum\limits_{V\in\mathfrak{F}} \sum\limits_{i,j=1}^{d_{[V]}}
d_{[V]}^{-1} \int_K  \Bigabs{\widehat{(g^V_{i
j})}(-\varphi_{\kappa}(\alpha))}^2\, d\lambda_K(\kappa)
.\end{equation}

Fix $\epsilon>0$, and set
$\norm{\mathfrak{F}}:=\sum_{V\in\mathfrak{F}} d_{[V]}$.  By the
Riemann--Lebesgue lemma for Abelian groups \cite[Proposition
4.13]{Fol}, there exists a symmetric compact subset $C$ of $\hA$
such that
\[ \Bigabs{\widehat{(g^V_{i j})}(\alpha)} <
\frac{\epsilon}{\norm{\mathfrak{F}}^{1/2}} \qquad\mbox{if }
\alpha\in\hA\smallsetminus C ,\] for all $V\in\mathfrak{F}$ and
$1\ls i,j\ls d_{[V]}$.  Then
\[ \Bigabs{\widehat{(g^V_{i
j})}(-\varphi_\kappa(\alpha))} <
\frac{\epsilon}{\norm{\mathfrak{F}}^{1/2}} \qquad\mbox{if}\
\alpha\in\hA\smallsetminus\hC ,\] for all $V\in\mathfrak{F}$,
$1\ls i,j\ls d_{[V]}$, and all $\kappa\in K$, where $\hC :=
\bigcup_{\kappa\in K} \varphi_\kappa(C)$; furthermore, since $C$
is compact so is $\hC$, by the continuity of the mapping
$(\alpha,\kappa)\mapsto\varphi_\kappa(\alpha)$. Assertion \1\ now
follows from \eqref{app-3}.

\smallskip

\noindent \2\  Fix $\epsilon>0$ and $\alpha\in\hA$. First, it
suffices to only consider functions $g$ of the form \eqref{app-1}
again. Fix such a function $g$ and write
$\mathcal{E}_\mathfrak{F}$ for the finite-dimensional subspace of
$L^2(K)$ spanned by the functions $V_{i j},\
i,j=1,\ldots,d_{[V]},\ V\in\mathfrak{F}$. We then claim that
$\varLambda_\alpha(g)(\mathcal{E}_\mathfrak{F}^\perp)=\set{0}$.
Indeed, inserting \eqref{app-2.5} into \eqref{app-1.5} yields that
\[[\varLambda_\alpha(g)\phi](\kappa) =
\sum\limits_{V\in\mathfrak{F}}\sum\limits_{i,j=1}^{d_{[V]}}
\widehat{(g^V_{i j})}(-\varphi_\kappa(\alpha))\cdot (V_{i
j}\ast\phi)(\kappa)\qquad (\kappa\in K)\] for any $\phi\in
L^2(K)$, whence $[\varLambda_\alpha(g)\phi]=0$ when $\phi$ is a
representative function of any irreducible unitary representation
of $K$ not equivalent to a representation in $\mathfrak{F}$, by
the Shur orthogonality relations.  Suppose then that
$\varLambda_\alpha = \bigoplus_{i\in\mathcal{I}} U_i$ with each
$U_i$ an irreducible unitary representation of $G$. Then
$\varLambda_\alpha(g) = \bigoplus_{i\in\mathcal{I}} U_i(g)$, and
since $\mathcal{E}_\mathfrak{F}$ is finite-dimensional, one must
have that $U_i(g)=0$ for all but finitely many $i$.
\end{proof}

\begin{proof}[Proof of Theorem \textup{\ref{app-thm-2}}]
If $\alpha,\beta\in\hA$ are characters of $A$, and $a\in A$, then
\[ \abs{\langle a,\alpha \rangle - \langle a,\beta \rangle}^2 =
\abs{1 - \langle a,\beta - \alpha \rangle}^2 .\] Thus, if $\mu\in
M(G)$ and $\phi\in L^2(K)$, then
\begin{align*}
\bigabs{[\hmu(\varLambda_\alpha-\varLambda_\beta) \phi]
(\tkappa)}^2 &= \biggabs{\int_{A\times K} [\langle
-a,\varphi_{\kappa\tkappa}(\alpha) \rangle - \langle
-a,\varphi_{\kappa\tkappa}(\beta) \rangle]\, \phi(\kappa\tkappa)\,
d\mu(a,\kappa)}^2\\ &\ls \norm{\mu} \int_{A\times K} \abs{1 -
\langle a,\varphi_{\kappa\tkappa}(\beta - \alpha) \rangle}^2
\abs{\phi(\kappa\tkappa)}^2\, d \lvert\mu\rvert(a,\kappa) ,
\end{align*}
for all $\tkappa\in K$, whence
\begin{align}\label{app-5}
 &\norm{\hmu(\varLambda_\alpha -
\varLambda_\beta) \phi}_{L^2(K)}^2  \\ &\qquad\ls \norm{\mu}
\int_K \int_{A\times K} \abs{1 - \langle
a,\varphi_{\kappa\tkappa}(\beta - \alpha) \rangle}^2
\abs{\phi(\kappa\tkappa)}^2\, d \lvert\mu\rvert(a,\kappa)\,
d\lambda_K(\tkappa) \nonumber \\ &\qquad = \norm{\mu} \int_K
\int_{A\times K} \abs{1 - \langle a,\varphi_{\tkappa}(\beta -
\alpha) \rangle}^2 \abs{\phi(\tkappa)}^2\, d\lambda_K(\tkappa)\, d
\abs{\mu}(a,\kappa) \nonumber \\ &\qquad = \norm{\mu} \int_K
\left[ \int_{A\times K} \abs{1 - \langle a,\varphi_{\tkappa}(\beta
- \alpha) \rangle}^2 d \lvert\mu\rvert(a,\kappa)\right]\,
\abs{\phi(\tkappa)}^2\, d\lambda_K(\tkappa)\nonumber
\end{align}
for $\alpha , \beta \in \hA$ and $\phi\in L^2(K)$, the first
equality by the left-invariance of Haar measure $\lambda_K$.

Next choose a compact set $A_1\subseteq A$ such that
$\lvert\mu\rvert ((A_1\times K)^c)<\epsilon$. Define
$\varPhi\colon A\times K\to A\times K$ by
$\varPhi(a,\kappa):=(\varphi_\kappa^{-1}(a),\kappa)$, which by the
continuity of the mapping $(a,\kappa)\mapsto
\varphi_{\kappa^{-1}}(a)$ is continuous, and set $A_2:=\pi_A
(\varPhi (A_1\times K))$ and $A_3:=\bigcap_{\kappa\in K}
\varphi_\kappa (A_2)$, where $\pi_A$ is the projection
$\pi_A\colon A\times K\to A$. Observe that $A_3\supseteq A_1$,
whence
\[ \lvert\mu\rvert ((A_3\times K)^c)<\epsilon ,\]
and that $A_2$ is compact.  Finally, let
\[ V_\epsilon := \bigset{\alpha\in\hA\colon \abs{1-\langle a,\alpha \rangle }
 <\epsilon\ \forall\,a\in A_2} ;\]
by the compactness of $A_2$, $V_\epsilon$ is an open neighborhood
of $0\in\hA$ \cite[\S 1.2.6]{rudin2}. For $a\in A_3$ and $\alpha -
\beta \in V_\epsilon$, one then has that
\[ \abs{1-\langle a, \varphi_{\tkappa}(\alpha - \beta)\rangle}
=  \bigabs{1-\langle \varphi_{\tkappa}^{-1}(a) , \alpha -
\beta\rangle } < \epsilon\qquad \mbox{for all }\tkappa\in K ,\]
whence
\[ \int_{A\times K}
\abs{1 - \langle a,\varphi_{\tkappa}(\beta - \alpha) \rangle}^2\,
d \lvert\mu\rvert(a,\kappa) \ls \epsilon^2\,
\lvert\mu\rvert(A_3\times K) + \lvert\mu\rvert((A_3\times K)^c)
\ls \epsilon^2\norm{\mu} + \epsilon \] for all $\tkappa\in K$ when
$\alpha-\beta\in V_\epsilon$, which by \eqref{app-5} completes the
proof.
\end{proof}

\end{document}